\documentclass[12pt,british,english]{article}
\usepackage{tgpagella}
\usepackage[T1]{fontenc}
\usepackage[latin9]{inputenc}
\usepackage{geometry}
\usepackage{units}
\usepackage{amsmath}
\usepackage{amsthm}
\usepackage{amssymb}
\usepackage{esint}

\makeatletter
\newcommand{\lyxaddress}[1]{
\par {\raggedright #1
\vspace{1.4em}
\noindent\par}
}
\theoremstyle{plain}
\newtheorem{thm}{\protect\theoremname}[section]
  \theoremstyle{definition}
  \newtheorem{defn}[thm]{\protect\definitionname}

\@ifundefined{date}{}{\date{}}
\makeatother

\usepackage{babel}
  \addto\captionsbritish{\renewcommand{\definitionname}{Definition}}
  \addto\captionsbritish{\renewcommand{\theoremname}{Theorem}}
  \addto\captionsenglish{\renewcommand{\definitionname}{Definition}}
  \addto\captionsenglish{\renewcommand{\theoremname}{Theorem}}
  \providecommand{\definitionname}{Definition}
\providecommand{\theoremname}{Theorem}

\begin{document}

\title{\textbf{\Large{}Kodaira Families and Newton-Cartan Structures with
Torsion}}

\author{James Gundry}
\maketitle

\lyxaddress{\begin{center}
{\small{}Department of Applied Mathematics and Theoretical Physics,
Centre for Mathematical Sciences, Wilberforce Road, Cambridge, CB3
0WA}
\par\end{center}}
\begin{abstract}
We describe the induced geometry on several classes of Kodaira moduli
spaces of rational curves in twistor spaces. By constructing connections
and frames on the moduli spaces we build and review twistor theories
pertaining to relativistic and non-relativistic geometries. Focussing
on the cases of three- and five-dimensional moduli spaces we establish
novel twistor theories of Newton-Cartan spacetimes. We generalise
a canonical class of connections on Kodaira moduli spaces to encompass
torsion and prove that in three dimensions deformations of the twistor
space's holomorphic structure induce Newton-Cartan structures with
torsion. In five dimensions the canonical connections contain generalised
Coriolis forces, and here also we consider the introduction of torsion
by deformation. Non-relativistic limits are exhibited as jumping phenomena
of normal bundles. We consider jumping phenomena in Newtonian twistor
theory, as well as exhibiting via twistor theory that three-dimensional
torsional Newton-Cartan geometry generically exists on the jumping
hypersurfaces of Gibbons-Hawking manifolds. \thispagestyle{empty}
\end{abstract}
\vfill{}

\begin{center}
{\small{}J.M.GUNDRY@DAMTP.CAM.AC.UK}
\par\end{center}{\small \par}

\textcompwordmark{}
\begin{center}
{\small{}DAMTP-2017-15}
\par\end{center}{\small \par}

\pagebreak{}

\tableofcontents{} \thispagestyle{empty}

\pagebreak{}

\section{Introduction}

\setcounter{page}{1}Let $Z$ be a complex manifold. One can study
the moduli space $M$ of a family of compact complex submanifolds
$X_{x}\subset Z$ for $x\in M$, and remarkably one finds that $M$
often comes equipped with canonical geometrical structures such as
metrics, forms, and affine connections.
\begin{defn}
\emph{\cite{MerkulovAppendices} A holomorphic family of compact complex
submanifolds embedded in a complex manifold $Z$ with moduli space
$M$ is a submanifold $F\subset Z\times M$ such that the restriction
$\nu:F\rightarrow M$ of the natural projection $Z\times M\rightarrow M$
is a proper regular map, and the submanifolds in $Z$ are $X_{x}=\mu\left(\nu^{-1}\left(x\right)\right)$
for $x\in M$, where $\mu:F\rightarrow Z$ is the restriction to $F$
of the natural projection $Z\times M\rightarrow Z$.}
\end{defn}
The submanifold $F$ is referred to as the correspondence space for
the construction. Kodaira showed when such a situation can occur,
giving a cohomological test for the submanifolds.
\begin{thm}
\cite{Kodaira,MerkulovAppendices} \label{thm:KodairaTheorem}Let
$X_{0}\subset Z$ be a compact complex embedded submanifold of a complex
manifold $Z$ such that $\check{H}^{1}\left(X_{x},N_{x}\right)=0$
(where $N_{x}=\nicefrac{TZ|_{X_{0}}}{TX_{0}}$ is its normal bundle).
Then $X_{0}$ is a member of a maximal holomorphic family $F$ with
moduli space $M$, and this family is complete in the sense that there
is an isomorphism $T_{x}M=\check{H}^{0}\left(X_{x},N_{x}\right).$
\end{thm}
The dimension of $M$ is thus equal to the dimension of $\check{H}^{0}\left(X_{x},N_{x}\right)$.
When the submanifolds are rational curves $X_{x}=\mathbb{P}^{1}$
we'll call the construction a twistor theory, and $Z$ is a twistor
space. The correspondence space is then $F=M\times\mathbb{P}^{1}$,
which can be identified with the projective primed spin bundle $P\mathbb{S}^{\prime}$,
the quotient by a homogeneity operator $\Upsilon=\pi_{A^{\prime}}\frac{\partial}{\partial\pi_{A^{\prime}}}$
of a rank-two holomorphic vector bundle $\mathbb{S}^{\prime}$ with
coordinates $\pi_{A^{\prime}}$ on fibres. The most famous example
of this construction is the nonlinear graviton construction of Penrose
\cite{PenroseNLGC}, in which $Z$ is a complex three-fold containing
a rational curve $X_{0}$ with normal bundle $N_{0}={\cal O}(1)\oplus{\cal O}(1)$;
$M$ is then a four-dimensional anti-self-dual conformal manifold,
and can be promoted to an Einstein manifold by imposing extra requirements.
The moduli space is naturally complex, but real slices with Euclidean
or neutral signature can be taken \cite{MDbook}. Other examples are
the cases
\begin{itemize}
\item $\dim Z=2$ with $N_{0}={\cal O}(2)$, where $M$ is a three-dimensional
Einstein-Weyl manifold \cite{HitchinEE};
\item $\dim Z=2$ with $N_{0}={\cal O}(1)$, where $M$ is a two-dimensional
projective manifold \cite{HitchinEE};
\item and $\dim Z=3$ with $N_{0}={\cal O}\oplus{\cal O}(2)$, where $M$
is a four-dimensional Newton-Cartan manifold \cite{DG16}.
\end{itemize}
In this paper we will add some new members to the list above by exhibiting
some novel Kodaira families with $X_{x}=\mathbb{P}^{1}$. Highlights
will include new Newtonian twistor spaces, and in section \ref{subsec:Galilean-Structures-and}
we'll prove the following.

\textcompwordmark{}

\textbf{Theorem \ref{thm:conformalNCin3d}.}\textbf{\emph{ }}\emph{Let
$Z={\cal O}\oplus{\cal O}(1)$ with global sections $X_{x}$. The
moduli space $M\ni x$ of these rational curves is a complex three-dimensional
manifold equipped with a family of Newton-Cartan structures parametrised
by three arbitrary functions on $M$ and an element of $\text{GL}\left(2,\mathbb{C}\right)$.}

\emph{\textcompwordmark{}}

We thus establish a flat model for the twistor theory of three-dimensional
Newton-Cartan spacetimes. A certain class of deformations of the complex
structure turn out to be interesting, inducing torsion. In section
\ref{subsec:Deformations-and-Torsion} we will prove the following.

\textbf{Theorem} \textbf{\ref{thm:2+1Torsion}.}\emph{ Let $Z$ be
the total space of an affine line bundle on ${\cal O}(1)$ with trivial
underlying translation bundle whose patching is 
\[
\hat{T}=T+f
\]
where $f$ represents a cohomology class in $\check{H}^{1}\left({\cal O}(1),{\cal O}_{{\cal O}(1)}\right)$.
The three-parameter family of global sections $X_{x}$ have normal
bundle $X_{x}={\cal O}\oplus{\cal O}(1)$ and the moduli space $M$
of those sections is a complex three-dimensional manifold equipped
with a family of torsional Newton-Cartan structures parametrised by
two arbitrary one-forms on $M$, two functions on $M$, and an element
of $\text{GL}\left(2,\mathbb{C}\right)$.}

\textcompwordmark{}

To prove this one must generalise the construction of affine connections
on Kodaira moduli spaces to include torsion; we make this generalisation
in section \ref{subsec:The-Torsion--Connection}. Torsion in Newton-Cartan
geometry is a subject which has attracted significant recent interest
\cite{Bergshoeff,BekMor,NoetherTNC}, and is \foreignlanguage{british}{characterised}
by the non-closure of the Newton-Cartan clock (see section \ref{subsec:Newton-Cartan-Geometry}
for an introduction to the relevant definitions). Analogous theorems
for the case of five-dimensional Newton-Cartan spacetimes will also
be proven, for which the twistor spaces are four-dimensional and the
isomorphism class of the normal bundles is ${\cal O}\oplus{\cal O}(1)\oplus{\cal O}(1)$.

We'll also discuss some Kodaira families whose geometry is not Newton-Cartan,
considering the induced geometry on families of rational curves with
normal bundles ${\cal O}(2n)$, ${\cal O}(2n-1)$, and ${\cal O}(2n-1)\oplus{\cal O}(2n-1)$
for $n\geq1$. The emphasis will be on the construction of the frame
on $M$, as well as two canonical families of affine connections.
Some of these examples have been considered previously in the literature,
as will be described in section \ref{sec:Some-Simple-Kodaira}, and
are included as examples of the frame and connection constructions
advocated in this paper. We will also discuss some of the links between
the above examples in terms of jumping phenomena.

The layout of the paper is as follows. In the remainder of this section
we will describe the methods of constructing induced geometry on Kodaira
families which will be preferred in this paper: we will consider the
construction of frame one-forms, and following Merkulov \cite{MerkulovAppendices},
the construction of two kinds of canonical affine connections. We
then proceed to generalise one of these, the $\Xi$-connection, to
include torsion. In section \ref{sec:Some-Simple-Kodaira} we'll apply
these ideas to some simple Kodaira families, finding ``relativistic''
induced geometry in various numbers of dimensions and discussing some
Newtonian limits. Section \ref{sec:NTT3} will describe the development
of a twistor theory of three-dimensional Newton-Cartan manifolds,
including the discussion of the introduction of torsion by Kodaira
deformations. Section \ref{sec:Novel4D} will introduce a couple of
supplementary results on Newtonian twistor theory in four dimensions.
Finally in section \ref{sec:Five-Dimensions} the material will be
extended to five dimensions, where we will find that the induced Newton-Cartan
geometry comes naturally equipped with (self-dual) Coriolis forces.

\subsection{Induced Metrics\label{subsec:Induced-Metrics}}

Given a Kodaira family of rational curves there are various ways of
constructing the induced geometry on $M$, all relying on the isomorphism
from theorem \ref{thm:KodairaTheorem}. In Penrose's original approach
one looks at the explicit alpha surfaces induced on $M$ by the pull-back
of the twistor coordinates to the correspondence space and then uses
the alpha surfaces to single out preferred vectors \cite{PenroseNLGC}.
Alternatively, in the case of the nonlinear graviton, one can directly
find the contravariant conformal structure as (the image on $M$ of)
the zero section of $\check{H}^{0}\left(Z,TZ\odot TZ\right)$. Yet
a third approach would be to construct an integral formula from twistor
cohomology classes and then compute the differential operator(s) on
$M$ for which the twistor classes constitute the kernel.

In this paper we'll employ yet another method, in which (families
of) frames of one-forms are directly computed on $M$.
\begin{thm}
\label{thm:vielbeintheorem}Let $M$ be the moduli space of a complete
holomorphic family of rational curves $X_{x}=\mathbb{P}^{1}$ in a
complex manifold $Z\rightarrow\mathbb{P}^{1}$ with normal bundles
$N_{x}\rightarrow\mathbb{P}^{1}$. $M$ is equipped with a preferred
family of one-forms - the frame - defined uniquely up to an element
of $\check{H}^{0}\left(X_{x},N_{x}\otimes N_{x}^{*}\right)$ per section
$X_{x}$.
\end{thm}

\subsubsection*{Proof}

Let $Z$ be $\left(k+1\right)$-dimensional. The family of one-forms
on the moduli space $M$ arises as a section of $N_{x}\otimes\Lambda_{x}^{1}\left(M\right)$
for each $x\in M$. Cover $\mathbb{P}^{1}$ by two patches $U$ and
$\hat{U}$ with coordinates $\lambda$ and $\hat{\lambda}$ respectively,
with holomorphic transition function 
\[
\hat{\lambda}=\lambda^{-1}
\]
on $U\cap\hat{U}$. We can describe $Z\overset{\pi}{\rightarrow}\mathbb{P}^{1}$
concretely as a complex manifold by exhibiting its patching. If $\hat{w}^{\mu}\left(w^{\nu},\lambda\right)$
is this patching (for $\mu,\nu=0,1,...,k$) then we can describe the
global section $X_{x}$ (over $U$) by the equation $w^{\mu}=w^{\mu}|\left(x,\lambda\right)$
for $x\in M$ parametrising the space of sections and where $w^{\mu}|$
are functions extracted from the patching.

One then finds that
\[
d\hat{w}^{\mu}|={\cal F}_{\,\nu}^{\mu}dw^{\nu}|
\]
where 
\[
{\cal F}_{\,\nu}^{\mu}\left(x,\lambda\right)=\frac{\partial\hat{w}^{\mu}}{\partial w^{\nu}}\left(w^{\alpha}|,\lambda\right)
\]
is the patching of the normal bundle $N_{x}$. (For notational convenience
we will use a vertical slash to denote the restriction to $X_{x}$
throughout this paper.) By the Birkhoff-Grothendieck theorem \cite{BG}
we can write
\[
{\cal F}=\hat{H}\,\text{diag}\left(\lambda^{-n_{1}},\lambda^{-n_{2}},...\,,\lambda^{-n_{k}}\right)H^{-1}
\]
where $H$ and $\hat{H}$ are holomorphic maps to $\text{GL}\left(k,\mathbb{C}\right)$
from $U$ and $\hat{U}$ respectively, and where the integers $(n_{1},n_{2},...,n_{k})$
specify the isomorphism class of $N$. We then extract the frame from
a section $\left(H^{-1}\right)_{\nu}^{\mu}dw^{\nu}|$ (or equivalently
$\left(\hat{H}^{-1}\right)_{\nu}^{\mu}d\hat{w}^{\nu}|$), and the
redundancy consists in multiplying $H$ (and $\hat{H}$) by a global
section of $N_{x}\otimes N_{x}^{*}$ (which may vary arbitrarily with
$x\in M$). We denote
\[
v\left(x\right)=H^{-1}dw|\in\check{H}^{0}\left(F|_{x},N_{x}\otimes\Lambda_{x}^{1}\left(M\right)\right)
\]
the \emph{frame section.}

The frame section then gives rise to a collection of one-forms $e^{AB...CA^{\prime}B^{\prime}...C^{\prime}}$
on $M$, where the unprimed indices arise from the component-structure
of $v$ and $A^{\prime}B^{\prime}C^{\prime}$ (each running from $0^{\prime}$
to $1^{\prime}$) arise from the dependence of $v$ on the base $\mathbb{P}^{1}$.
For example, in the case where 
\[
N_{x}={\cal O}(2)\oplus{\cal O}(2)
\]
we write
\[
v^{A}=e^{AA_{1}^{\prime}A_{2}^{\prime}}\pi_{A_{1}^{\prime}}\pi_{A_{2}^{\prime}}\,\,,
\]
where $\left[\pi_{A^{\prime}}\right]$ are homogeneous coordinates
on the base $\mathbb{P}^{1}$ and where $A=0,1$. One then extracts
the frame $e^{AA_{1}^{\prime}A_{2}^{\prime}}=e_{a}^{AA_{1}^{\prime}A_{2}^{\prime}}\left(x^{b}\right)dx^{a}$.\hfill{}$\square$

In this paper we will include in $v$ the result of the most general
element of $\check{H}^{0}\left(X_{x},N_{x}\otimes N_{x}^{*}\right)$
per line $X_{x}$ in the guise of arbitrary functions. (This is analogous
to writing a conformal structure $\left[g\right]$ as a single metric
$g=\alpha g_{0}$ for some representative $g_{0}\in\left[g\right]$
and some arbitrary non-vanishing function $\alpha$. In some cases
this analogy is an identity.)
\begin{defn}
\emph{All tensor fields on $M$ which can be constructed from the
frame using only the tensor product and the symplectic forms $\epsilon_{AB}$
and $\epsilon_{A^{\prime}B^{\prime}}$ are induced on $M$ as the
span of $v$.\label{def:spanofv}}
\end{defn}
Often, pleasingly, the span of $v$ contains a metric (which may depend
on some number of arbitrary functions), and in the nonlinear graviton
construction this metric is exactly the conformal structure we could
induce in a more traditional way (i.e. via the presence of alpha surfaces
in $M$ as is done in \cite{PenroseNLGC}).

Using alpha surfaces to find the induced conformal structure (or otherwise)
becomes increasingly complicated in higher dimensions, but the frame
method of theorem \ref{thm:vielbeintheorem} does not. We will thus
make frequent use of it, though without losing sight completely of
the existence of alpha surfaces. 

\subsection{Induced Affine Connections\label{subsec:Induced-Affine-Connections}}

The construction of affine connections on the moduli spaces of complete
holomorphic families of submanifolds was considered by Merkulov \cite{MerkGeom,MerkulovAppendices}.
This work led to the solution of the holonomy problem \cite{MerkSchwHol}.

First consider torsion-free affine connections more generally. Let
$J_{x}^{\,k}$ be the ideal of germs of holomorphic functions on $M$
which vanish to order $k$ at $x\in M$. The second-order tangent
bundle $T^{[2]}M$ is defined to be the union over all points in $M$
of second-order tangent spaces
\[
T_{x}^{[2]}M=\left(\nicefrac{J_{x}}{J_{x}^{\,3}}\right)^{*}.
\]
An element $\left(V^{ab},V^{a}\right)$ of $T_{x}^{[2]}M$ consists
of the first two non-vanishing terms of the Taylor expansion of a
function vanishing at $x$; a section of $T^{[2]}M$ gives rise to
a second-order linear differential operator
\[
V^{[2]}=\left(V^{ab},V^{a}\right)\rightsquigarrow V^{ab}\partial_{a}\partial_{b}+V^{a}\partial_{a}\,\,,
\]
where for brevity's sake we put $\partial_{a}=\frac{\partial}{\partial x^{a}}$.
There is a short exact sequence
\begin{equation}
0\rightarrow TM\rightarrow T^{[2]}M\rightarrow\odot^{2}TM\rightarrow0\label{TMsequence}
\end{equation}
with maps
\[
\left(V^{a}\right)\mapsto\left(0,V^{a}\right)\qquad\text{and}\qquad\left(V^{ab},V^{a}\right)\mapsto\left(V^{ab}\right).
\]
A torsion-free affine connection $\nabla$ on $TM$ is then equivalent
to a (left) splitting of (\ref{TMsequence}), i.e. a linear map
\begin{equation}
\gamma\,:\,T^{[2]}M\rightarrow TM\label{connectionmaps}
\end{equation}
acting as
\[
\gamma\,:\,\left(V^{ab},V^{a}\right)\mapsto\left(V^{a}+\Gamma_{\,bc}^{a}V^{bc}\right)
\]
for functions $\Gamma_{\,bc}^{a}$ on $M$ which we identify as Christoffel
symbols. We'll now describe, following \cite{MerkulovAppendices},
two ways of constructing such a map in twistor theory. The original
treatment in \cite{MerkulovAppendices} is considerably more sophisticated
than what is given here, as Merkulov describes the construction for
a general Kodaira moduli space. In the case $Z\rightarrow\mathbb{P}^{1}$
the construction is much simpler.

\subsubsection{The $\Xi$-Connection\label{subsec:The-Xi-Connection}}

As above cover $\mathbb{P}^{1}$ with two open sets $U$ and $\hat{U}$
with respective coordinate functions $\lambda$ and $\hat{\lambda}$
subject to the transition function $\hat{\lambda}=\lambda^{-1}$ on
$U\cap\hat{U}$. These coordinates are known as \emph{inhomogeneous}
coordinates on $\mathbb{P}^{1}$. Again consider a general fibred
twistor space $\pi:Z\rightarrow\mathbb{P}^{1}$ characterised by the
holomorphic patching
\begin{equation}
\hat{w}^{\mu}=\hat{w}^{\mu}\left(w^{\nu},\lambda\right)\,\,,\label{eq:patching}
\end{equation}
where $w^{\mu}$ are the twistor coordinates on the fibres. Henceforth
restrict to cases in which one has 
\[
\check{H}^{1}\left(\mathbb{P}^{1},N_{x}\otimes N_{x}^{*}\right)=0.
\]

A section $V=\left(V^{ab},V^{a}\right)$ of $T^{\left[2\right]}M$
gives rise to a differential operator $V^{ab}\partial_{a}\partial_{b}+V^{a}\partial_{a}$
which following \cite{MerkulovAppendices} we'll now apply to (\ref{eq:patching}).
\begin{equation}
V^{ab}\partial_{a}\partial_{b}\hat{w}^{\mu}+V^{a}\partial_{a}\hat{w}^{\mu}=V^{ab}{\cal F}_{\,\nu}^{\mu}\partial_{a}\partial_{b}w^{\nu}+V^{a}{\cal F}_{\,\nu}^{\mu}\partial_{a}w^{\nu}+V^{ab}{\cal F}_{\,\nu\rho}^{\mu}\partial_{a}w^{\nu}\partial_{b}w^{\rho},\label{equationfour}
\end{equation}
where again we have
\[
{\cal F}_{\,\nu}^{\mu}=\frac{\partial\hat{w}^{\mu}}{\partial w^{\nu}}|
\]
and
\[
{\cal F}_{\,\nu\rho}^{\mu}=\frac{\partial^{\mu}\hat{w}}{\partial w^{\nu}\partial w^{\rho}}|.
\]
(Recall that a vertical slash indicates that one must restrict to
global sections $w^{\nu}=w^{\nu}|\left(x^{a},\lambda\right)$.) If
we can write this line as a global section of $N_{x}$ then we have
(via the Kodaira isomorphism $T_{x}M=\check{H}^{0}\left(\mathbb{P}^{1},N_{x}\right)$)
constructed a map $\gamma$. The only problematic term is the last
one. The $\Xi$-connection is constructed by splitting 
\begin{equation}
{\cal F}_{\,\nu\rho}^{\mu}\partial_{a}w^{\nu}=-\hat{\chi}_{\,\alpha\,a}^{\mu}{\cal F}_{\,\rho}^{\alpha}+{\cal F}_{\,\alpha}^{\mu}\chi_{\,\rho\,a}^{\alpha}\label{eq: Xi-splitting}
\end{equation}
for a 0-cochain $\left\{ \chi\right\} $ of $N_{x}\otimes N_{x}^{*}\otimes T_{x}^{*}M$
per point $x\in M$. The left-hand-side of (\ref{eq: Xi-splitting})
is always (in the $Z\rightarrow\mathbb{P}^{1}$ case) a 1-cocycle
of $N_{x}\otimes N_{x}^{*}\otimes T_{x}^{*}M$, and since by assumption
$\check{H}^{1}\left(\mathbb{P}^{1},N_{x}\otimes N_{x}^{*}\right)=0$
the splitting is always possible.

Equation (\ref{equationfour}) then becomes
\[
V^{ab}\partial_{a}\partial_{b}\hat{w}^{\mu}+V^{a}\partial_{a}\hat{w}^{\mu}+V^{ab}\hat{\chi}_{\,\alpha\,(a}^{\mu}\partial_{b)}\hat{w}^{\alpha}={\cal F}_{\,\nu}^{\mu}\left(V^{ab}\partial_{a}\partial_{b}w^{\nu}+V^{a}\partial_{a}w^{\nu}+V^{ab}\chi_{\,\rho\,(a}^{\nu}\partial_{b)}w^{\rho}\right)
\]
and so we have constructed a global section of $N_{x}$ per point
$x$. The connection symbols for the map thus constructed can be read
off as
\[
\partial_{a}\partial_{b}w^{\nu}+\chi_{\,\rho\,(a}^{\nu}\partial_{b)}w^{\rho}=\Gamma_{\,ab}^{c}\partial_{c}w^{\nu}.
\]
There is, though, a possible source of non-uniqueness. If $\check{H}^{0}\left(\mathbb{P}^{1},N_{x}\otimes N_{x}^{*}\right)\neq0$
then one is free to add any element of that group to both $\hat{\chi}_{\,\alpha\,a}^{\mu}$
and $\chi_{\,\alpha\,a}^{\mu}$. Therefore what one obtains is an
equivalence class of connections (the $\Xi$-connection). As with
the frame we will choose to describe the equivalence class by a single
most-general representative containing arbitrary functions.

\subsubsection*{Example}

Consider $Z={\cal O}(1)\oplus{\cal O}(1)$ with patching
\[
\hat{\Omega}^{A}=\lambda^{-1}\Omega^{A}\qquad\hat{\lambda}=\lambda
\]
for $A=0,1$ and global sections
\[
\Omega^{A}|=x^{A0^{\prime}}\lambda+x^{A1^{\prime}}
\]
where $x^{AA^{\prime}}$ are coordinates on $M$ parametrising the
family of global sections. (This is the twistor space for flat four-dimensional
spacetime.)

The Christoffel symbols for the $\Xi$-connection are then given by
\begin{equation}
\Gamma_{\qquad CC^{\prime}DD^{\prime}}^{AA^{\prime}}=\chi_{D\,\,CC^{\prime}}^{A}\delta_{D^{\prime}}^{A^{\prime}}+\chi_{C\,\,DD^{\prime}}^{A}\delta_{C^{\prime}}^{A^{\prime}}\,\,\,,\label{eq:XiforO1O1}
\end{equation}
where $\chi_{\,B}^{A}$ are four arbitrary one-forms on $M$ constituting
a global section of $N_{x}\otimes N_{x}^{*}\otimes T_{x}^{*}M$ per
point $x\in M$.

\subsubsection{The $\Lambda$-Connection}

An alternative twistor construction of a class of maps (\ref{connectionmaps})
gives rise to the so-called $\Lambda$-connection, which we now briefly
\foreignlanguage{british}{summarise}. This class of connections often
degenerates into a single affine connection, and in as much as it
is necessary and useful to make such a distinction, it is the $\Lambda$-connection
which should be considered the \emph{physical }connection.

Consider again the patching for a general fibred twistor space $Z\rightarrow\mathbb{P}^{1}$:
\[
\hat{w}^{\mu}=\hat{w}^{\mu}\left(w^{\nu},\lambda\right)\,,
\]
where $w^{\mu}$ are the coordinates on the fibres, and again consider
the equation (\ref{equationfour}) resulting from the action of the
section $V$ of $T^{\left[2\right]}M$. To construct the $\Lambda$-connection
we do the splitting differently.

We instead choose to solve
\begin{equation}
{\cal F}_{\,\alpha\beta}^{\mu}=-\hat{\sigma}_{\,\nu\rho}^{\mu}{\cal F}_{\,\alpha}^{\nu}{\cal F}_{\,\beta}^{\rho}+{\cal F}_{\,\eta}^{\mu}\sigma_{\,\alpha\beta}^{\eta}\label{splitting problem}
\end{equation}
for a $0$-cochain $\left\{ \sigma\right\} $ of $N\otimes\left(N^{*}\odot N^{*}\right)\rightarrow\mathbb{P}^{1}$,
and the Christoffel symbols for the resulting map $\gamma$ can be
read off from
\begin{equation}
\Gamma_{\,bc}^{a}\partial_{a}w^{\mu}|=\partial_{b}\partial_{c}w^{\mu}|+\sigma_{\,\nu\rho}^{\mu}\partial_{b}w^{\nu}|\partial_{c}w^{\rho}|\label{readoffchristoffel}
\end{equation}
(or the equivalent expression over $\hat{U}$). In the $Z\rightarrow\mathbb{P}^{1}$
case the left hand side of (\ref{splitting problem}) is always a
1-cocycle of $N\otimes\left(N^{*}\odot N^{*}\right)$. 

The difficult part of this construction is the solution of the splitting
problem (\ref{splitting problem}), which in some cases is not possible,
and is often not unique.

Uniqueness is determined by whether there are global sections of $N\otimes\left(N^{*}\odot N^{*}\right)$;
if these exist then one is free to add one to $\left\{ \sigma\right\} $
and so construct a different connection. In Penrose's case we have
\[
\check{H}^{0}\left(\mathbb{P}^{1},N\otimes\left(N^{*}\odot N^{*}\right)\right)=0
\]
and so the connection is always unique. This is unsurprising; we can
always call upon the Levi-Civita connection.

There are Kodaira deformations (giving rise to ${\cal F}_{\,\alpha\beta}^{\mu}$)
for which (\ref{splitting problem}) cannot be solved iff
\[
\check{H}^{1}\left(\mathbb{P}^{1},N\otimes\left(N^{*}\odot N^{*}\right)\right)\neq0\,,
\]
and there are several reasons this may occur. One is that the spacetime
suffers a jump in the normal bundle; another is when the torsion-free
requirement essential to the construction is broken. In Penrose's
case we can calculate that $\check{H}^{1}\left(\mathbb{P}^{1},N\otimes\left(N^{*}\odot N^{*}\right)\right)$
vanishes, so all Kodaira deformations lead to torsion-free connections,
in line with the nonlinear graviton construction.

\subsubsection{The Torsion $\Xi$-Connection\label{subsec:The-Torsion--Connection}}

In sections \ref{subsec:Deformations-and-Torsion} and \ref{subsec:FiveDDeformations}
the $\Lambda$-connections fail to exist and the $\Xi$-connections
cannot be made compatible with the induced frame data. The frame section
suggests that the reason for this is that the moduli space's connection
possesses torsion. In this subsection we will generalise the $\Xi$-connection
of \cite{MerkulovAppendices,MerkGeom} to include torsion.

Consider a general (possibly torsional) affine connection to be a
map
\[
\gamma:\,\Gamma\left(TM\times TM\right)\rightarrow\Gamma\left(TM\right).
\]
As in the constructions of the previous two sections, we'll build
such a map from the twistor data via the Kodaira isomorphism. Let
$V=V^{a}\partial_{a}$ and $W=W^{a}\partial_{a}$ be vector fields
on $M$. We want to use the complex structure of the twistor space
$Z$ to build a torsional connection by directly constructing $\nabla_{W}V=\gamma\left(V,W\right)$.

Apply $V$ to the patching (\ref{eq:patching}) to obtain

\begin{equation}
V^{a}\partial_{a}\hat{w}^{\mu}={\cal F}_{\,\nu}^{\mu}V^{a}\partial_{a}w^{\nu}\label{eq:Kexample}
\end{equation}
 as usual. Then apply $W$ to (\ref{eq:Kexample}) to obtain
\[
W^{b}\partial_{b}V^{a}\partial_{a}\hat{w}^{\mu}+W^{b}V^{a}\partial_{b}\partial_{a}\hat{w}^{\mu}=W^{b}\partial_{b}{\cal F}_{\,\nu}^{\mu}V^{a}\partial_{a}w^{\nu}+W^{b}{\cal F}_{\,\nu}^{\mu}\partial_{b}V^{a}\partial_{a}w^{\nu}+W^{b}{\cal F}_{\,\nu}^{\mu}V^{a}\partial_{b}\partial_{a}w^{\nu}.
\]
As in the torsion-free case the only obstruction to this line (evaluated
at $x\in M$) constituting a global section of $N_{x}$ is the first
term on the right-hand-side, and one must decide what to do with it.

If we can write
\begin{equation}
\partial_{b}{\cal F}_{\,\nu}^{\mu}=-\hat{\rho}_{\,\alpha\,b}^{\mu}{\cal F}_{\,\nu}^{\alpha}+{\cal F}_{\,\beta}^{\mu}\rho_{\,\nu\,b}^{\beta}\label{eq:this}
\end{equation}
for some $0$-cochain $\left\{ \rho\right\} $ of $N_{x}\otimes\Lambda_{x}^{1}\left(M\right)$
for each $x\in M$ then
\[
W^{b}\partial_{b}V^{a}\partial_{a}\hat{w}^{\mu}+W^{b}V^{a}\partial_{b}\partial_{a}\hat{w}^{\mu}+W^{b}\hat{\rho}_{\,\alpha\,b}^{\mu}V^{a}\partial_{a}\hat{w}^{\alpha}={\cal F}_{\,\beta}^{\mu}\left(W^{b}\rho_{\,\nu\,b}^{\beta}V^{a}\partial_{a}w^{\nu}+W^{b}\partial_{b}V^{a}\partial_{a}w^{\beta}+W^{b}V^{a}\partial_{b}\partial_{a}w^{\beta}\right)
\]
constitutes a global section of $N_{x}$ (for each $x\in M$), and
hence a vector field on $M$ via the Kodaira isomorphism. One can
then extract the connection symbols from
\begin{equation}
\Gamma_{\,ab}^{c}\partial_{c}w^{\mu}|=\partial_{a}\partial_{b}w^{\mu}|+\rho_{\,\nu\,b}^{\mu}\partial_{a}w^{\nu}|\label{eq:torsionconnectionreadoff}
\end{equation}
(or its counterpart over $\hat{U}$) just as in the torsion-free case.

The connection symbols arising from (\ref{eq:torsionconnectionreadoff})
generically possess torsion, and the torsion-free part of the connection
agrees with the torsion-free $\Xi$-connection of exhibited in section
\ref{subsec:The-Xi-Connection}. We accordingly call the connection
of this section the \emph{torsion $\Xi$-connection. }

Just like the $\Xi$-connection the existence is determined by the
non-vanishing of $\check{H}^{1}\left(X_{x},N_{x}\otimes N_{x}^{*}\right)$
and the connection is defined up to an element of $\check{H}^{0}\left(X_{x},N_{x}\otimes N_{x}^{*}\right)\otimes\Lambda_{x}^{1}\left(M\right)$
per $x\in M$, giving us a family of connections induced on $M$.

\subsection{Newton-Cartan Geometry\label{subsec:Newton-Cartan-Geometry}}

The non-relativistic limit of a Lorentzian manifold is a Newton-Cartan
manifold; such manifolds provide a geometrical setting for non-relativistic
physics \cite{Cartan}. As in general relativity a four-manifold models
the spacetime, and test particles travel on geodesics of a (usually
torsion-free) connection. There's a metric also, though unlike in
general relativity the metric and the connection are independent.
In this section we will describe Newton-Cartan manifolds in some detail,
taking \cite{DuvalHorvathy} as a reference.

Before introducing the full Newton-Cartan geometry we'll begin with
a subordinate definition.
\begin{defn}
\emph{A $\left(d+1\right)$-dimensional Galilean spacetime is a triple
$(M,\,h,\,\theta)$ where}
\end{defn}
\begin{itemize}
\item \emph{$M$ is a $\left(d+1\right)$-dimensional manifold;}
\item \emph{$h$ is a symmetric tensor field of valence $\begin{pmatrix}2\\
0
\end{pmatrix}$ with signature $(0++...+)$ (and so has rank $d$) called the metric;}
\item \emph{and $\theta$ is a closed one-form spanning the kernel of $h$
called the clock.}
\end{itemize}
The pair $(h,\theta)$ is called a \emph{Galilean structure, }and
the number of spatial dimensions is $d$.

Since $\theta$ is closed we can always locally write $\theta=dt$
for some function $t:M\rightarrow\mathbb{R}$. This function is then
taken as a coordinate on the time axis, a one-dimensional submanifold
over which the spacetime $M$ is fibred. We call the \foreignlanguage{british}{fibres}
\emph{spatial slices} and when restricted to such a slice the metric
$h$ is a more familiar signature $(+...++)$ $d$-metric. Throughout
this paper the indices $a,b,c$ will run from $0$ to $d$ and the
spatial indices $i,j,k$ will run from $1$ to $d$.
\begin{defn}
\emph{A $\left(d+1\right)$-dimensional Newton-Cartan spacetime is
a quadruple $(M,\,h,\,\theta,\,\nabla)$ where}
\end{defn}
\begin{itemize}
\item \emph{$M$ is a $(d+1)$-dimensional manifold;}
\item \emph{$(h,\theta)$ is a Galilean structure; }
\item \emph{and $\nabla$ is a torsion-free connection compatible with the
Galilean structure in the sense that $\nabla h=0$ and $\nabla\theta=0$.}
\end{itemize}
Note that crucially $\nabla$ must be specified independently of the
metric and clock. 

The field equations for Newton-Cartan gravity arise as the Newtonian
limit of the Einstein equations \cite{Limits} and are given by
\begin{equation}
R_{ab}=4\pi G\rho\theta_{a}\theta_{b}\label{eq:NCfieldequations}
\end{equation}
where $R_{ab}$ is the Ricci tensor associated to $\nabla$; $G$
is Newton's constant; and $\rho:M\rightarrow\mathbb{R}$ is the mass
density. In addition to the field equations there is the Trautman
condition \cite{DuvalHorvathy}
\begin{equation}
h^{a[b}R_{\,\,(de)a}^{c]}=0\,\,,\label{TrautmanCondition}
\end{equation}
where $R_{\,\,bcd}^{a}$ is the Riemann tensor of $\nabla$ which
ensures that there exist potentials (such as the Newtonian potential)
for the connection components. Newton-Cartan connections which satisfy
(\ref{TrautmanCondition}) are called \emph{Newtonian} connections.

The field equations imply that $h$ is flat on spatial slices, and
we can always introduce \emph{Galilean} coordinates $(t,x^{i})$ such
that
\begin{equation}
h=\delta^{ij}\frac{\partial}{\partial x^{i}}\otimes\frac{\partial}{\partial x^{j}}\qquad\mbox{and}\qquad\theta=dt\label{standardGalileanstructure}
\end{equation}
for $i=1,2,...,d$. We'll refer to (\ref{standardGalileanstructure})
as the \emph{standard} Galilean structure.

Only connections compatible with $\theta$ and $h$ are allowed by
definition; one can show \cite{Dautcourt} that the most general such
connection has components
\begin{equation}
\Gamma_{\,\,bc}^{a}=\frac{1}{2}h^{ad}\left(\partial_{b}h_{cd}+\partial_{c}h_{bd}-\partial_{d}h_{bc}\right)+\partial_{(b}\theta_{c)}U^{a}+\theta_{(b}F_{c)d}h^{ad}\label{eq:generalNCconnection}
\end{equation}
where
\begin{itemize}
\item $U^{a}$ is any vector field satisfying $\theta(U)=1$;
\item $F_{ab}$ is any two-form;
\item and $h_{ab}$ is uniquely determined by $h^{ab}h_{bc}=\delta_{c}^{a}-\theta_{c}U^{a}$
and $h_{ab}U^{b}=0$.
\end{itemize}
Possible connections, given a Galilean structure, are then \foreignlanguage{british}{determined}
by a choice of $(U,\,F)$. The Trautman condition (\ref{TrautmanCondition})
is equivalent to the statement that $F$ is closed, and hence for
a Newtonian connection we can locally write $F=dA$. As well as the
obvious gauge symmetry
\[
A\,\longmapsto\,A+d\chi
\]
there is a further redundancy in this description. There exist \emph{Milne
boosts }which can be thought of as gauge transformations of $(U,\,F)$
leaving $\Gamma_{\,\,bc}^{a}$ unchanged \cite{DuvalHorvathy}.

With $d=3$ the most general vacuum Newton-Cartan manifold satisfying
(\ref{eq:NCfieldequations}) and (\ref{TrautmanCondition}) then has
\[
\Gamma_{\,\,tt}^{i}=\delta^{ij}\partial_{j}V\qquad\mbox{and}\qquad\Gamma_{\,\,jt}^{i}=\Gamma_{\,\,tj}^{i}=\delta_{jl}\epsilon^{ilk}\partial_{k}\varOmega
\]
\begin{equation}
\mbox{where}\qquad\delta^{ij}\partial_{i}\partial_{j}V+2\delta^{ij}\partial_{i}\varOmega\partial_{j}\varOmega=0\qquad\mbox{and}\qquad\delta^{ij}\partial_{i}\partial_{j}\varOmega=0,\label{eq:NCfieldeqns}
\end{equation}
with all other connection components vanishing. The corresponding
two-form $F$ is given by
\[
F=-dV\wedge dt\,+\,\epsilon_{ijk}\delta^{kl}\partial_{l}\varOmega\,dx^{i}\wedge dx^{j}.
\]
The geodesic equations suggest interpreting the function $V$ as the
Newtonian (gravitational) potential and the function $\varOmega$
as a potential for \foreignlanguage{british}{generalised} (spatially-varying)
Coriolis forces. Note that although the degrees of freedom in a Newton-Cartan
connection appear similar to that of an electromagnetic field the
equations (\ref{eq:NCfieldeqns}) governing them are more complicated.

Of recent interest has been \emph{torsional }Newton-Cartan geometry,
in which the connection is allowed to have some torsion. This is manifest
in the clock failing to be closed because $d\theta\neq0$ is incompatible
with $\nabla\theta=0$ for a torsion-free connection. Equation (\ref{eq:generalNCconnection})
is modified to include the skew part of $\partial_{a}\theta_{b}$
too, giving rise to the torsion. Later in this paper we will construct
Kodaira families with induced torsional Newton-Cartan structures,
featuring clocks which aren't closed.

\pagebreak{}

\section{Examples of Frames and Connections\label{sec:Some-Simple-Kodaira}}

In this section we'll apply the results of subsections \ref{subsec:Induced-Metrics}
and \ref{subsec:Induced-Affine-Connections} in some simple cases.

\subsection{Line Bundles on $\mathbb{P}^{1}$}

Twistor spaces with families of submanifolds $X_{x}$ having $N_{x}={\cal O}(k)$
for some $k\geq1$ are about as simple as it gets; they are, though,
sufficiently sophisticated as to require some careful treatment of
their canonical connections, particularly when $k$ is odd.

Applied in these cases theorem \ref{thm:vielbeintheorem} amounts
to the construction of a \emph{paraconformal} structure on $M,$ i.e.
a bundle isomorphism 
\[
TM=\odot^{k}\mathbb{S^{\prime}}
\]
as is studied in \cite{DunajskiTodParaconformal,BaileyEastwoodPara,Bryant},
concretely given by the frame $e_{a}^{A_{1}^{\prime}...A_{k}^{\prime}}$.

\subsubsection{Odd Dimensions}

When $k$ is even the treatment is relatively straightforward: the
span of the frame will give us a conformal structure and the $\Lambda$-connection
can pick out a preferred representative.
\begin{thm}
Let $Z\rightarrow\mathbb{P}^{1}$ be a complex two-fold containing
a rational curve $X_{0}$ with normal bundle $N_{0}={\cal O}\left(2n\right)$
for some $n\geq1$. The Kodaira moduli space of rational curves $X_{x}$
is a $\left(2n+1\right)$-dimensional complex conformal manifold.\label{thm:O(2n)}
\end{thm}
For $n=1$ $Z$ is a minitwistor space \cite{HitchinEE}.

\subsubsection*{Proof}

Since $\check{H}^{1}\left(\mathbb{P}^{1},{\cal O}\left(2n\right)\right)=0$
the rational curve $X_{0}$ is a member of a $\dim\check{H}^{0}\left(\mathbb{P}^{1},{\cal O}\left(2n\right)\right)=\left(2n+1\right)$-dimensional
family of rational curves, and by theorem \ref{thm:vielbeintheorem}
we obtain a section of $\Lambda_{x}^{1}\left(M\right)\otimes N_{x}$
at each point $x\in M$ which gives rise to a frame via
\[
v=e^{A_{1}^{\prime}...A_{2n}^{\prime}}\pi_{A_{1}^{\prime}}...\pi_{A_{2n}^{\prime}}
\]
 and so a metric
\begin{equation}
g=e_{A_{1}^{\prime}...A_{2n}^{\prime}}\otimes e^{A_{1}^{\prime}...A_{2n}^{\prime}}\label{eq:metricO(2n)}
\end{equation}
of maximal rank in the span of $v$. The redundancy acts as
\[
e^{A_{1}^{\prime}...A_{2n}^{\prime}}\mapsto\alpha e^{A_{1}^{\prime}...A_{2n}^{\prime}}
\]
for any $\alpha:M\rightarrow\mathbb{C}^{*}$, resulting in conformal
transformations $g\mapsto\alpha^{2}g$.\hfill{}$\square$

In the case for which the patching for $Z$ is that of ${\cal O}(2n)$
(even if the sections are deformed) one may fix a particular metric
from the conformal class by constructing the $\Lambda$-connection,
which is in this case unique and exists for the ${\cal O}(2n)$ patching.

We can equip $Z$ with an involution which singles out Euclidean signature
metrics. (See, for example, \cite{DGT16,MDbook}.) The metric (\ref{eq:metricO(2n)})
is the same as that arising from the classical invariant theory described
in \cite{DunajskiPenrose}.

\subsubsection{Even Dimensions}

When $k$ is odd the situation is more complicated because the span
contains no (non-degenerate) metric. In the following theorem we consider
one option of what one \emph{can} do with the frame, though this is
by no means the only geometry induced on $M$.
\begin{thm}
Let $Z\rightarrow\mathbb{P}^{1}$ be a complex two-fold containing
a rational curve $X_{0}$ with normal bundle $N_{0}={\cal O}\left(2n-1\right)$
for some $n\geq1$. Then the Kodaira moduli space $M$ of rational
curves $X_{x}$ is a $\left(2n\right)$-dimensional complex torsional
projective manifold.\label{thm:O(2n-1)}
\end{thm}
Restricting to torsion-free connections only, for $n=1$ this is the
standard twistor theory of projective surfaces due to Hitchin \cite{HitchinEE}
and for $n=2$ the normal bundle ${\cal O}(3)$ is that associated
to exotic holonomies in the work of Bryant \cite{Bryant}, whose twistor
theory is described in terms of solutions spaces of ODEs in \cite{DunajskiTodParaconformal}.

\subsubsection*{Proof}

Since $\check{H}^{1}\left(\mathbb{P}^{1},{\cal O}\left(2n-1\right)\right)=0$
(for $n\geq1$) the rational curve $X_{0}$ is a member of a $\dim\check{H}^{0}\left(\mathbb{P}^{1},{\cal O}\left(2n-1\right)\right)=\left(2n\right)$-dimensional
family of rational curves, and by theorem \ref{thm:vielbeintheorem}
we obtain a section of $\Lambda_{x}^{1}\left(M\right)\otimes N_{x}$
at each point $x\in M$ which gives rise to a frame via
\[
v=e^{A_{1}^{\prime}...A_{2n-1}^{\prime}}\pi_{A_{1}^{\prime}}...\pi_{A_{2n-1}^{\prime}}.
\]
Unlike the case of odd dimensions the span does not contain a metric
of maximal rank. We can, though, construct a family of connections
out of the frame. A change of global section of $N_{x}\otimes N_{x}^{*}$
acts as 
\begin{equation}
v\mapsto\alpha v\label{eq:gauge}
\end{equation}
 for $\alpha:M\rightarrow\mathbb{C}^{*}$, so write the frame as
\[
e^{A_{1}^{\prime}...A_{2n-1}^{\prime}}=\alpha\left(x\right)\varsigma_{a}^{A_{1}^{\prime}...A_{2n-1}^{\prime}}\left(x\right)dx^{a}.
\]
We can construct a canonical family of affine connections on $M$
by requiring $\nabla e^{A_{1}^{\prime}...A_{2n-1}^{\prime}}=0$. Concretely,
this gives us
\begin{equation}
\Gamma_{\,ab}^{c}=\varsigma_{A_{1}^{\prime}...A_{2n-1}^{\prime}}^{c}\partial_{a}\varsigma_{b}^{A_{1}^{\prime}...A_{2n-1}^{\prime}}+\delta_{a}^{c}\partial_{b}\ln\alpha\label{eq:oddkconnectionsymbols}
\end{equation}
where $\varsigma_{A_{1}^{\prime}...A_{2n-1}^{\prime}}^{a}$ is the
inverse of $\varsigma_{a}^{A_{1}^{\prime}...A_{2n-1}^{\prime}}$.
The connections described in (\ref{eq:oddkconnectionsymbols}) possess
torsion whenever $de^{A_{1}^{\prime}...A_{2n-1}^{\prime}}\neq0$;
their torsion-free parts (and hence their geodesics) constitute a
projective structure, in that a change of $\alpha$ leaves the unparametrised
geodesics unaltered.

\hfill{}$\square$

Now consider the canonical connections induced on $M$ without reference
to the frame section. We have $N_{x}\otimes N_{x}^{*}={\cal O}$,
so the torsion $\Xi$-connection always exists and depends on a single
one-form on $M$. In the case $Z={\cal O}(1)$ the torsion-free $\Xi$-connection
is a standard flat projective structure.

On the other hand we have
\[
N_{x}\otimes\left(N_{x}^{*}\odot N_{x}^{*}\right)={\cal O}(1-2n)
\]
so
\[
\check{H}^{0}\left(\mathbb{P}^{1},N_{x}\otimes\left(N_{x}^{*}\odot N_{x}^{*}\right)\right)=0
\]
and
\[
\check{H}^{1}\left(\mathbb{P}^{1},N_{x}\otimes\left(N_{x}^{*}\odot N_{x}^{*}\right)\right)=\mathbb{C}^{2n-2}.
\]
Thus the $\Lambda$-connection, when it exists, is unique.

For $Z={\cal O}(1)$ we find that $\Gamma_{\,bc}^{a}=0$, so the moduli
space comes equipped with a preferred representative of the projective
structure, and moreover one which is metrisable. There is thus in
this case an important corollary: $M$ is equipped with a flat metric
$h_{ab}$. We simply impose $\nabla h=0$ and the torsion-free condition
(by analogy with the existence of the Levi-Civita connection), giving
us a metric with constant coefficients (which is unique up to diffeomorphisms
in two dimensions). This will be important for theorem \ref{thm:conformalNCin3d}.
(Note that this does not imply that all such $\Lambda$-connections
give rise to metrics: the connection is not guaranteed to be metrisable.)

In theorem \ref{thm:O(2n-1)} we chose to make the whole frame parallel,
but we had other options. Another would be to construct a family of
connections by declaring the form $e_{A_{1}^{\prime}...A_{2n-1}^{\prime}}\otimes e^{A_{1}^{\prime}...A_{2n-1}^{\prime}}$
to be parallel. In two dimensions this form is complex symplectic,
and the connection is known as a symplectic connection.

\subsection{${\cal O}\left(2n-1\right)\oplus{\cal O}\left(2n-1\right)\rightarrow\mathbb{P}^{1}$}
\begin{thm}
Let $Z\rightarrow\mathbb{P}^{1}$ be a complex three-fold containing
a rational curve $X_{0}$ with normal bundle $N_{0}={\cal O}\left(2n-1\right)\oplus{\cal O}\left(2n-1\right)$
for some $n\geq1$. Then the Kodaira moduli space of rational curves
$X_{x}$ is a $\left(4n\right)$-dimensional complexified conformal
manifold.\label{thm:4n}
\end{thm}
For $n=1$ $Z$ is a standard twistor space, and in a different context
this class of normal bundles is the setting for the heavenly hierarchy
described in \cite{DunajskiMasonhK}.

Note that theorem \ref{thm:4n} is a different construction of $4n$-dimensional
moduli spaces to that in \cite{Salamon,PedersenPoon}, where the authors
induce quaternionic structures on Kodaira families of global sections
of manifolds with normal bundle $\oplus^{2n}{\cal O}\left(1\right)$.

\subsubsection*{Proof}

Since 
\[
\check{H}^{1}\left(\mathbb{P}^{1},{\cal O}\left(2n-1\right)\oplus{\cal O}\left(2n-1\right)\right)=0
\]
 the rational curve $X_{0}$ is a member of a family of dimension
\[
\dim\check{H}^{0}\left(\mathbb{P}^{1},{\cal O}\left(2n-1\right)\oplus{\cal O}\left(2n-1\right)\right)=4n\,\,,
\]
and by theorem \ref{thm:vielbeintheorem} we obtain a section of $\Lambda_{x}^{1}\left(M\right)\otimes N_{x}$
at each point $x\in M$ which gives rise to a frame $e^{AA_{1}^{\prime}...A_{2n-1}^{\prime}}$
via
\[
v^{A}=e^{AA_{1}^{\prime}...A_{2n-1}^{\prime}}\pi_{A_{1}^{\prime}}...\pi_{A_{2n-1}^{\prime}}
\]
 and so a metric
\[
g=e_{AA_{1}^{\prime}...A_{2n-1}^{\prime}}\otimes e^{AA_{1}^{\prime}...A_{2n-1}^{\prime}}.
\]
The redundancy acts via an invertible global section of $N_{x}\otimes N_{x}^{*}$,
which takes
\[
\begin{pmatrix}e^{00^{\prime}...0^{\prime}}\lambda^{2n-1}+e^{00^{\prime}...0^{\prime}1^{\prime}}\lambda^{2n-2}+...+e^{01^{\prime}...1^{\prime}}\\
e^{10^{\prime}...0^{\prime}}\lambda^{2n-1}+e^{10^{\prime}...0^{\prime}1^{\prime}}\lambda^{2n-2}+...+e^{11^{\prime}...1^{\prime}}
\end{pmatrix}\mapsto\begin{pmatrix}a_{0} & b_{0}\\
c_{0} & d_{0}
\end{pmatrix}\begin{pmatrix}e^{00^{\prime}...0^{\prime}}\lambda^{2n-1}+e^{00^{\prime}...0^{\prime}1^{\prime}}\lambda^{2n-2}+...+e^{01^{\prime}...1^{\prime}}\\
e^{10^{\prime}...0^{\prime}}\lambda^{2n-1}+e^{10^{\prime}...0^{\prime}1^{\prime}}\lambda^{2n-2}+...+e^{11^{\prime}...1^{\prime}}
\end{pmatrix}
\]
for any four functions $\phi_{\,B}^{A}=\begin{pmatrix}a_{0} & b_{0}\\
c_{0} & d_{0}
\end{pmatrix}:\mathbb{C}\rightarrow\text{GL}\left(2,\mathbb{C}\right)$, resulting in 
\[
e^{AA_{1}^{\prime}...A_{2n-1}^{\prime}}\mapsto\phi_{\,B}^{A}e^{BA_{1}^{\prime}...A_{2n-1}^{\prime}}
\]
and so
\[
g\mapsto\epsilon_{AD}\phi_{\,B}^{D}\phi_{\,C}^{A}e_{\,\,A_{1}^{\prime}...A_{2n-1}^{\prime}}^{B}\otimes e^{CA_{1}^{\prime}...A_{2n-1}^{\prime}}=\left(\det\phi\right)g.
\]
We thus obtain a conformal transformation.\hfill{}$\square$

Again the $\Lambda$-connection, when it exists, can be used to fix
a particular representative $g\in\left[g\right]$ by imposing $\nabla g=0$.
(The $\Lambda$-connection exists in this case when the deformation
giving rise to $Z$ does not fall within $\check{H}^{1}\left(X_{x},N_{x}\otimes\left(N_{x}^{*}\odot N_{x}^{*}\right)\right)$
when restricted to twistor lines.)

\subsection{Limits in $4n$ Dimensions}

Let $D=4n$ where $n>0$ is an integer.
\begin{thm}
Let $Z_{c}\rightarrow\mathbb{P}^{1}$ be a one-parameter family of
vector bundles with patching
\[
\hat{T}=T+\frac{S}{c\lambda^{2n-1}}
\]
\[
\hat{S}=\lambda^{2-D}S.
\]
For $c\neq\infty$ the normal bundle to all rational curves $X_{x}$
is $N_{x}={\cal O}\left(2n-1\right)\oplus{\cal O}\left(2n-1\right)$
and the homogeneous frame section is
\[
v=\begin{pmatrix}v^{0}\\
v^{1}
\end{pmatrix}\qquad\text{where}\qquad v^{A}=e^{AA_{1}^{\prime}...A_{2n-1}^{\prime}}\pi_{A_{1}^{\prime}...A_{2n-1}^{\prime}}
\]
giving rise to a non-degenerate metric
\[
g=e_{AA_{1}^{\prime}...A_{2n-1}^{\prime}}\otimes e^{AA_{1}^{\prime}...A_{2n-1}^{\prime}}
\]
on $M$. For $c=\infty$ the normal bundle to all rational curves
$X_{x}$ is $N_{x}={\cal O}\oplus{\cal O}\left(4n-2\right)$ and the
homogeneous frame section is
\[
v=\begin{pmatrix}\theta\\
e^{A_{1}^{\prime}...A_{4n-2}^{\prime}}\pi_{A_{1}^{\prime}}...\pi_{A_{4n-2}^{\prime}}
\end{pmatrix}
\]
giving rise to a Galilean structure with clock $\theta$ and
\[
h^{-1}=e_{A_{1}^{\prime}...A_{4n-2}^{\prime}}\otimes e^{A_{1}^{\prime}...A_{4n-2}^{\prime}}.
\]
\end{thm}
The induced geometry is subject to a redundancy, which in the $c\neq\infty$
case amounts to a conformal ambiguity and in the $c=\infty$ constitutes
the non-metric nature of the connection's gravitational sector. For
$n=1$ this is the standard Newtonian limit of twistor theory presented
in \cite{DG16}, and for $c=\infty$ the manifold is a Newton-Cartan
manifold with arbitrary gravitational sector.

\subsubsection*{Proof}

We need to begin by identifying the isomorphism class of $N_{x}$,
which will be the same for all $x\in M$ because $Z_{c}$ is the total
space of a vector bundle. For $c=\infty$ the patching for $N_{x}$
is 
\[
{\cal F}_{\,}=\begin{pmatrix}1 & 0\\
0 & \lambda^{2-D}
\end{pmatrix}
\]
and so the isomorphism class is obvious: $N_{x}={\cal O}\oplus{\cal O}\left(D-2\right)$.
We thus obtain a frame section
\[
v=H^{-1}\begin{pmatrix}dT|\\
dS|
\end{pmatrix}
\]
where $H^{-1}$ is a general global section of $N_{x}\otimes N_{x}^{*}$.
The rational curves are given by
\[
T|=t\qquad S|=x_{0}+x_{1}\lambda+...+x_{D-2}\lambda^{D-2}
\]
and span of $v$ is the Galilean structure which was advertised.

For $c\neq\infty$ the patching of $N_{x}$ splits as
\[
{\cal F}_{\,}=\hat{H}\begin{pmatrix}\lambda^{1-2n} & 0\\
0 & \lambda^{1-2n}
\end{pmatrix}H^{-1}
\]
where (for instance) 
\[
H=\begin{pmatrix}1 & 0\\
-c\lambda^{2n-1} & 1
\end{pmatrix}\qquad\text{and}\qquad\hat{H}=\begin{pmatrix}0 & c^{-1}\\
-c & \hat{\lambda}^{2n-1}
\end{pmatrix}.
\]
This exhibits the normal bundle's isomorphism class as ${\cal O}\left(2n-1\right)\oplus{\cal O}\left(2n-1\right)$,
and the frame section is
\[
v=H_{0}\begin{pmatrix}1 & 0\\
c\lambda^{2n-1} & 1
\end{pmatrix}\begin{pmatrix}dT|\\
dS|
\end{pmatrix}=H_{0}\begin{pmatrix}dT|\\
dS|+c\lambda^{2n-1}dT|
\end{pmatrix}
\]
for an arbitrary non-degenerate matrix of functions $H_{0}$. The
rational curves are given by
\[
T|=t-\frac{1}{c}\sum_{m=0}^{2-2n}x_{m+2n}\lambda^{1+m}
\]
\[
S|=x_{0}+x_{1}\lambda+...+x_{D-2}\lambda^{D-2}\,,
\]
which results in a frame section of the advertised form.\hfill{}$\square$

\pagebreak{}

\section{Newtonian Twistor Theory in Three Dimensions\label{sec:NTT3}}

The twistor theory of complexified three-dimensional manifolds with
non-degenerate connections is called \emph{minitwistor} theory and
is well understood \cite{HitchinEE}. In this section we will consider
the twistor theory of three-dimensional Newton-Cartan manifolds. The
relevant twistor spaces are three-dimensional and will be \foreignlanguage{british}{characterised}
by the normal bundle to twistor lines $X_{x}$ being
\[
N_{x}={\cal O}\oplus{\cal O}(1).
\]

Families of rational curves with normal bundles isomorphic to ${\cal O}\oplus{\cal O}(1)$
have been considered in \cite{Gindikin}, where it is described that
the tangent spaces of the three-dimensional moduli space come equipped
with preferred one-parameter families of null rays. As we shall see
below, this makes the isomorphism class ${\cal O}\oplus{\cal O}(1)$
well-suited to describing Newton-Cartan structures in three dimensions.

It is straightforward to see that
\[
\check{H}^{1}\left(\mathbb{P}^{1},N_{x}\right)=0\qquad\text{and}\qquad\check{H}^{0}\left(\mathbb{P}^{1},N_{x}\right)=\mathbb{C}^{3}
\]
and so a three-dimensional moduli space is feasible. Additionally
\[
\check{H}^{1}\left(\mathbb{P}^{1},N_{x}\otimes N_{x}^{*}\right)=0\qquad\text{and}\qquad\check{H}^{0}\left(\mathbb{P}^{1},N_{x}\otimes N_{x}^{*}\right)=\mathbb{C}^{4}
\]
so that the isomorphism class is stable; the $\Xi$-connection always
exists and always depends on four arbitrary one-forms. Finally we
see that
\[
\check{H}^{1}\left(\mathbb{P}^{1},N_{x}\otimes\left(N_{x}^{*}\odot N_{x}^{*}\right)\right)=\mathbb{C}\qquad\text{and}\qquad\check{H}^{0}\left(\mathbb{P}^{1},N_{x}\otimes\left(N_{x}^{*}\odot N_{x}^{*}\right)\right)=\mathbb{C}^{4}
\]
so that the $\Lambda$-connection fails to exist for some Kodaira
deformations, and when it does exist it is not unique, depending on
four arbitrary functions. The case in which the $\Lambda$-connection
fails to exist is when the deformation introduces torsion.

We'll begin by considering the undeformed case $Z={\cal O}\oplus{\cal O}(1)$.
After calculating the canonical connections and the Galilean structure
we will then proceed to deform $Z$.

\subsection{Galilean Structures and Canonical Connections\label{subsec:Galilean-Structures-and}}

In this subsection we'll discuss the canonical geometry induced on
the moduli space of global sections of $Z={\cal O}\oplus{\cal O}(1)$.
The patching is
\[
\hat{T}=T\qquad\hat{\Omega}=\lambda^{-1}\Omega
\]
where $\hat{\lambda}=\lambda^{-1}$ as usual for the base $\mathbb{P}^{1}$.
The global sections are
\[
w^{\mu}|=\begin{pmatrix}T|\\
\Omega|
\end{pmatrix}=\begin{pmatrix}t\\
y+z\lambda
\end{pmatrix}
\]
for $x^{a}=\left(t,y,z\right)\in\mathbb{C}^{3}=M$. (Recall that a
vertical slash indicates the restriction to rational curves.) We'll
also wish to use homogeneous coordinates $\left[\pi_{A^{\prime}}\right]$
on the base and $\omega$ for the ${\cal O}(1)$ fibre, writing
\[
\omega|=x^{A^{\prime}}\pi_{A^{\prime}}
\]
for the global sections.
\begin{thm}
Let $Z={\cal O}\oplus{\cal O}(1)$ with global sections $X_{x}$.
The moduli space $M\ni x$ of these rational curves is a complex three-dimensional
manifold equipped with a family of Newton-Cartan structures parametrised
by three arbitrary functions on $M$ and an element of $\text{GL}\left(2,\mathbb{C}\right)$.\label{thm:conformalNCin3d}
\end{thm}
Two of these functions determine the gravitational field; the other
is a conformal factor.

\subsubsection*{Proof}

To begin we will construct the frame induced on $M$ using theorem
\ref{thm:vielbeintheorem}, which is done by splitting the patching
${\cal F}$ for the normal bundle; we must solve
\[
\begin{pmatrix}1 & 0\\
0 & \lambda^{-1}
\end{pmatrix}H=\hat{H}\begin{pmatrix}1 & 0\\
0 & \lambda^{-1}
\end{pmatrix}
\]
for
\[
H=\begin{pmatrix}h_{1} & h_{2}\\
h_{3} & h_{4}
\end{pmatrix}\qquad\text{and}\qquad\hat{H}=\begin{pmatrix}\hat{h}_{1} & \hat{h}_{2}\\
\hat{h}_{3} & \hat{h}_{4}
\end{pmatrix}
\]
as holomorphic maps from $U$ and $\hat{U}$ to $\text{GL}(2,\mathbb{C})$.
This amounts to the four individual splitting problems
\[
\hat{h}_{1}=h_{1}\qquad\hat{h}_{2}=\lambda h_{2}
\]
\[
\hat{h}_{3}=\lambda^{-1}h_{3}\qquad\hat{h}_{4}=h_{4}
\]
whose general solution is
\[
H=\begin{pmatrix}m & 0\\
a_{0}+a_{1}\lambda & k
\end{pmatrix}
\]
for four arbitrary holomorphic functions $\left(m,k,a_{0},a_{1}\right)$
on $M$ constrained only by $m\neq0$ and $k\neq0$. The frame section
can then be read off from $v=H^{-1}dw|$, giving us a clock
\begin{equation}
\theta=m^{-1}dt\label{eq:2+1clock}
\end{equation}
and spatial one-forms
\[
e^{0^{\prime}}=k^{-1}\left(dx^{0^{\prime}}-m^{-1}a_{1}dt\right)
\]
\[
e^{1^{\prime}}=k^{-1}\left(dx^{1^{\prime}}-m^{-1}a_{0}dt\right).
\]

To proceed further we must calculate the $\Lambda$-connection on
$M$ as described in section \ref{subsec:Induced-Metrics}. The splitting
problem to be solved is
\[
0=-\hat{\sigma}_{\,\nu\rho}^{\mu}{\cal F}_{\,\alpha}^{\nu}{\cal F}_{\,\beta}^{\rho}+{\cal F}_{\,\nu}^{\mu}\sigma_{\,\alpha\beta}^{\nu}
\]
for a $0$-cochain $\left\{ \sigma\right\} $ valued in $N_{x}\otimes\left(N_{x}^{*}\odot N_{x}^{*}\right)$
on each twistor line $X_{x}$, which amounts to
\[
\hat{\sigma}_{\,TT}^{T}=\sigma_{\,TT}^{T}\qquad\hat{\sigma}_{\,T\Omega}^{T}=\lambda\sigma_{\,T\Omega}^{T}\qquad\hat{\sigma}_{\,\Omega\Omega}^{T}=\lambda^{2}\sigma_{\,\Omega\Omega}^{T}
\]
\[
\hat{\sigma}_{\,TT}^{\Omega}=\lambda^{-1}\sigma_{\,TT}^{\Omega}\qquad\hat{\sigma}_{\,T\Omega}^{\Omega}=\sigma_{\,T\Omega}^{\Omega}\qquad\hat{\sigma}_{\,\Omega\Omega}^{\Omega}=\lambda\sigma_{\,\Omega\Omega}^{\Omega}
\]
and so
\[
\sigma_{\,TT}^{T}=\Sigma\qquad\sigma_{\,TQ}^{T}=0
\]
\[
\sigma_{\,\Omega\Omega}^{T}=0\qquad\sigma_{\,TT}^{\Omega}=\phi_{0}+\lambda\phi_{1}
\]
\[
\sigma_{\,T\Omega}^{\Omega}=\chi\qquad\sigma_{\,\Omega\Omega}^{\Omega}=0
\]
for any four functions $\left(\Sigma,\chi,\phi_{0},\phi_{1}\right)$
on $M$. The connection symbols are then
\[
\Gamma_{\,tt}^{t}=\Sigma\qquad\Gamma_{\,it}^{t}=0\qquad\Gamma_{\,ij}^{t}=0
\]
\[
\Gamma_{\,tt}^{y}=\phi_{0}\qquad\Gamma_{\,tt}^{z}=\phi_{1}
\]
\[
\Gamma_{\,yt}^{y}=\chi\qquad\Gamma_{\,zt}^{y}=0\qquad\Gamma_{\,yt}^{z}=0\qquad\Gamma_{\,zt}^{z}=\chi
\]
\[
\Gamma_{\,jk}^{i}=0.
\]
The result is that the moduli space comes equipped with a family of
connections containing gravitational forces (described by the functions
$\phi_{0}$ and $\phi_{1}$).

The connection allows us to restrict the clock (\ref{eq:2+1clock})
by imposing $\nabla\theta=0$, which tells us that
\[
\Sigma=-\partial_{t}\ln m
\]
and
\[
\partial_{A^{\prime}}m=0\,\,,
\]
so that $m$ is a function of $t$ alone. Given that $m$ is now just
a non-vanishing function of $t$, we see that (\ref{eq:2+1clock})
is a standard Newton-Cartan clock, where the function $m$ just allows
for diffeomorphisms of the time axis, with $\Sigma$ ensuring that
upon such diffeomorphisms the clock remains parallel.

To complete the proof we must construct a family of Newton-Cartan
metrics. The data already induced defines a metric as follows, by
requiring the usual conditions $h\left(\theta,\,\,\right)=0$ and
$\nabla h=0$. These imply that we must have $h^{at}=0$ and that
$h^{ij}$ obeys
\[
\partial_{t}h^{ij}+2\chi h^{ij}=0
\]
and
\[
\partial_{k}h^{ij}=0.
\]
We deduce that $h^{ij}$ must be any element of $\text{GL}\left(2,\mathbb{C}\right)$
multiplied by an arbitrary non-vanishing function of $t$, and we
also have that $\chi=\chi\left(t\right)$ only. Constant non-degenerate
two-by-two metrics are all equal up to (restricted) diffeomorphisms
\[
y\mapsto\alpha y+\beta z\qquad z\mapsto\gamma y+\delta z
\]
for $\begin{pmatrix}\alpha & \beta\\
\gamma & \delta
\end{pmatrix}\in\text{GL}\left(2,\mathbb{C}\right)$, so we are free to take any such member $\tilde{h}^{ij}$ as our
metric, giving us
\[
h^{ij}=\kappa\left(t\right)\tilde{h}^{ij}
\]
where
\[
\kappa=\exp\left\{ -2\int\chi\,dt\right\} 
\]
is non-vanishing and determined by the arbitrary function $\chi$.

Thus we have a family of Galilean structures $\left(h,\theta\right)$
and a family of connections $\nabla$, two of whose arbitrary functions
describe gravitational forces, completing the proof.

\hfill{}$\square$

To fix a specific gravitational sector for the Newton-Cartan manifold
we require, as in \cite{DG16}, some additional data on $Z$. The
following theorem provides one way of specifying this data.
\begin{thm}
Equip $Z={\cal O}\oplus{\cal O}(1)$ with a $1$-cocycle taking values
in its canonical bundle $K\rightarrow Z$. This induces a preferred
global section of ${\cal O}(1)$ which can be used to fix a complex
Newton-Cartan structure out of the complex Newton-Cartan structure
of theorem \ref{thm:conformalNCin3d}, where the gravitational sector
is locally of the form $\Gamma_{\,tt}^{i}=g^{i}$ for $\boldsymbol{g}$
divergence-free and determined uniquely by $f$.
\end{thm}

\subsubsection*{Proof}

A simple calculation shows that $K={\cal O}(-3)_{Z}$, and so a $1$-cocycle
is represented by a function $f$ of weight minus three in homogeneous
coordinates, and provides a (Serre-dual) global section of ${\cal O}(1)$
by
\[
\phi_{A^{\prime}}\left(x\right)=\frac{1}{2\pi i}\oint_{\Gamma}\pi_{A^{\prime}}f\left(T|,\omega|,\pi_{B^{\prime}}\right)\pi\cdot d\pi
\]
where $\omega|=x^{A^{\prime}}\pi_{A^{\prime}}=y\pi_{1^{\prime}}+z\pi_{0^{\prime}}$.
We have
\[
\frac{\partial}{\partial x^{A^{\prime}}}\phi^{A^{\prime}}=0
\]
automatically. Taking this to fix $\sigma_{\,TT}^{\Omega}$ we have
\[
\Gamma_{\,tt}^{A^{\prime}}\partial_{A^{\prime}}\omega|=\phi^{A^{\prime}}\pi_{A^{\prime}}
\]
and so
\[
\Gamma_{\,tt}^{y}=\phi^{1^{\prime}}\qquad\Gamma_{\,tt}^{z}=\phi^{0^{\prime}}.
\]
Thus the gravitational sector is fixed to be a unique divergence-free
$\boldsymbol{g}$ given by the global $\phi^{A^{\prime}}$. We note
that the divergence-free condition ensures that the Newton-Cartan
spacetime is vacuum.

\hfill{}$\square$

\subsubsection*{Torsion-Free $\Xi$-Connection}

Theorem \ref{thm:conformalNCin3d} employed the $\Lambda$-connection
because it is the most powerful construction available in terms of
constraining moduli space geometry. It is interesting, however, to
consider the $\Xi$-connection also so as to compare the flat model
of theorem \ref{thm:conformalNCin3d} with the torsion-inducing deformations
of theorem \ref{thm:2+1Torsion} where the torsion $\Xi$-connection
is the only connection available.

Here we'll calculate the canonical torsion-free $\Xi$-connection
for $Z={\cal O}\oplus{\cal O}(1)$. The calculation is straightforward;
since $\partial_{a}{\cal F}_{\,\nu}^{\mu}=0$ one must solve
\[
0=-\hat{\chi}_{\,\nu\,\,a}^{\mu}{\cal F}_{\,\rho}^{\nu}+{\cal F}_{\,\nu}^{\mu}\chi_{\,\rho\,\,a}^{\nu}
\]
for the most general $0$-cochain $\left\{ \chi_{\,\nu\,a}^{\mu}\right\} $
of $N_{x}\otimes N_{x}^{*}\otimes\Lambda_{x}^{1}\left(M\right)$ for
each $x\in M$, where
\[
{\cal F}_{\,\nu}^{\mu}=\begin{pmatrix}1 & 0\\
0 & \lambda^{-1}
\end{pmatrix}.
\]
This becomes the four individual splitting problems given by
\[
\hat{\chi}_{\,T\,\,a}^{T}=\chi_{\,T\,\,a}^{T}\qquad\hat{\chi}_{\,\Omega\,\,a}^{T}=\lambda\chi_{\,\Omega\,\,a}^{T}
\]
\[
\hat{\chi}_{\,T\,\,a}^{\Omega}=\lambda^{-1}\chi_{\,T\,\,a}^{\Omega}\qquad\hat{\chi}_{\,\Omega\,\,a}^{\Omega}=\chi_{\,\Omega\,\,a}^{\Omega}\,\,,
\]
which have the general solution
\[
\chi_{\,T\,\,a}^{T}=C_{a}\qquad\chi_{\,\Omega\,\,a}^{T}=0
\]
\[
\chi_{\,T\,\,a}^{\Omega}=A_{a}^{0}+\lambda A_{a}^{1}\qquad\chi_{\,\Omega\,\,a}^{\Omega}=B_{a}
\]
for arbitrary one-forms $\left(A,B,C,\kappa\right)$ on $M$. We then
read-off the connection from
\[
\Gamma_{\,bc}^{a}\partial_{a}w^{\mu}|=\partial_{b}\partial_{c}w^{\mu}|+\chi_{\,\nu\,\,b}^{\mu}\partial_{c}w^{\nu}|+\chi_{\,\nu\,\,c}^{\mu}\partial_{b}w^{\nu}|
\]
which here leads to
\[
\Gamma_{\,tt}^{t}=2C_{t}\qquad\Gamma_{\,it}^{t}=C_{i}\qquad\Gamma_{\,ij}^{t}=0
\]
\[
\Gamma_{\,tt}^{y}=A_{t}^{0}\qquad\Gamma_{\,tt}^{z}=A_{t}^{1}
\]
\[
\Gamma_{\,yt}^{y}=A_{y}^{0}+B_{t}\qquad\Gamma_{\,zt}^{y}=A_{z}^{0}\qquad\Gamma_{\,zt}^{z}=A_{z}^{1}+B_{t}\qquad\Gamma_{\,yt}^{z}=A_{y}^{1}
\]
\[
\Gamma_{\,ij}^{k}=2B_{(i}\delta_{j)}^{k}.
\]
Thus the $\Xi$-connection comprises all connections which have $\Gamma_{\,ij}^{t}=0$
and have flat projective structures as their spatial sectors. Compatibility
with the (closed) clock $\theta=m^{-1}dt$ then imposes
\[
\Gamma_{\,ab}^{t}=\delta_{a}^{t}m\partial_{b}\left(m^{-1}\right)
\]
so (recalling that the connection is torsion-free) one must put $C_{i}=0$
and
\[
C_{t}=\frac{1}{2}m\partial_{t}\left(m^{-1}\right).
\]
The remaining freedom in $\left(\Gamma_{\,ab}^{c},\theta_{a}\right)$
is then given by three one-forms and one non-vanishing function on
the time axis.

\subsection{Deformations and Torsion\label{subsec:Deformations-and-Torsion}}

A natural next step is to consider deforming the complex structure
of $Z={\cal O}\oplus{\cal O}(1)$, and a case of interest is when
we write ${\cal O}\oplus{\cal O}(1)$ as a trivial affine line bundle
on ${\cal O}(1)$ and then deform the patching for the affine line
bundle so that we have
\begin{equation}
\hat{T}=T+f\left(\Omega,\lambda\right)\label{eq:2+1deformation}
\end{equation}
\[
\hat{\Omega}=\lambda^{-1}\Omega
\]
as the patching for $Z\rightarrow{\cal O}(1)\rightarrow\mathbb{P}^{1}$.
The analogous deformation in four-dimensional Newtonian twistor theory
\cite{DG16} leads to a jump in the isomorphism class of the normal
bundle to every\footnote{\emph{Almost} every twistor line; see section \ref{subsec:On-Jumping-Points}.}
twistor line from ${\cal O}\oplus{\cal O}(2)$ to ${\cal O}(1)\oplus{\cal O}(1)$.
In the three-dimensional case this cannot be what occurs, because
the isomorphism class of the normal bundle is stable.

The deformation leading to (\ref{eq:2+1deformation}), when restricted
to twistor lines, corresponds exactly to the part of $N_{x}\otimes\left(N_{x}^{*}\odot N_{x}^{*}\right)$
which causes $\check{H}^{1}\left(\mathbb{P}^{1},N_{x}\otimes\left(N_{x}^{*}\odot N_{x}^{*}\right)\right)$
to fail to vanish. Thus we expect something to go ``wrong'' with
the torsion-free affine connection on $M$; in fact what happens is
that the connection fails to be torsion-free.
\begin{thm}
\label{thm:2+1Torsion}Let $Z$ be the total space of an affine line
bundle on ${\cal O}(1)$ with trivial underlying translation bundle
whose patching is 
\[
\hat{T}=T+f
\]
where $f$ represents a cohomology class in $\check{H}^{1}\left({\cal O}(1),{\cal O}_{{\cal O}(1)}\right)$.
The three-parameter family of global sections $X_{x}$ have normal
bundle $X_{x}={\cal O}\oplus{\cal O}(1)$ and the moduli space $M$
of those sections is a complex three-dimensional manifold equipped
with a family of torsional Newton-Cartan structures parametrised by
two arbitrary one-forms on $M$, two functions on $M$, and an element
of $\text{GL}\left(2,\mathbb{C}\right)$.
\end{thm}

\subsubsection*{Proof}

The proof will proceed in stages.
\begin{enumerate}
\item Using theorem \ref{thm:vielbeintheorem} we will construct a clock
one-form $\theta$ on $M$ depending on one non-vanishing function
$m$ on $M$; this clock will not be closed, meaning that it cannot
be made compatible with a torsion-free connection.
\item We will then construct the torsion $\Xi$-connection described in
section \ref{subsec:The-Torsion--Connection}; its connection symbols
will depend on four arbitrary one-forms $\left(A^{0},A^{1},B,C\right)$
on $M$.
\item Imposing $\nabla\theta=0$ for the torsion $\Xi$-connection is then
possible, and results in the fixing of $C$.
\item The remaining piece of data, the Newton-Cartan metric $h$, will then
be constructed by imposing $h\left(\theta,\,\,\right)=0$ and $\nabla h=0$.
This restricts the remaining one-forms by a closure condition and
determines the metric up to a choice of constant two-by-two non-degenerate
matrix $\tilde{h}$.
\end{enumerate}
To begin we must construct the twistor functions by finding the global
sections of $Z\rightarrow\mathbb{P}^{1}$, as these are required for
theorem \ref{thm:vielbeintheorem}. It'll be useful to expand the
representative $f$ when restricted to sections of ${\cal O}(1)$.
For $\Omega|=y+z\lambda$ (with $y$ and $z$ coordinates on $M$)
we can write
\[
f\left(\Omega|,\lambda\right)=\sum_{n=-\infty}^{\infty}\gamma_{n}\lambda^{n}
\]
where $\gamma_{n}\left(y,z\right)$ are functions one can extract
via integration.

We then have
\[
\hat{T}=T+\sum_{n=-\infty}^{\infty}\gamma_{n}\lambda^{n}
\]
and so the twistor functions are
\[
T|=t-\sum_{n=1}^{\infty}\gamma_{n}\lambda^{n}\qquad\hat{T}|=t+\sum_{n=-\infty}^{0}\gamma_{n}\lambda^{n}
\]
where we have chosen to put the $\gamma_{0}$ term into $\hat{T}|$
(without loss of generality: we could always effect the diffeomorphism
$t\mapsto t-\gamma_{0}$). Stage one of the proof is then to use these
twistor functions to calculate a frame section via theorem \ref{thm:vielbeintheorem}.
This involves solving the splitting problem
\[
\begin{pmatrix}1 & \frac{\partial f}{\partial\Omega}|\\
0 & \lambda^{-1}
\end{pmatrix}\begin{pmatrix}h_{1} & h_{2}\\
h_{3} & h_{4}
\end{pmatrix}=\begin{pmatrix}\hat{h}_{1} & \hat{h}_{2}\\
\hat{h}_{3} & \hat{h}_{4}
\end{pmatrix}\begin{pmatrix}1 & 0\\
0 & \lambda^{-1}
\end{pmatrix}
\]
where as usual 
\[
H=\begin{pmatrix}h_{1} & h_{2}\\
h_{3} & h_{4}
\end{pmatrix}\qquad\text{and}\qquad\hat{H}=\begin{pmatrix}\hat{h}_{1} & \hat{h}_{2}\\
\hat{h}_{3} & \hat{h}_{4}
\end{pmatrix}
\]
constitute holomorphic maps to $\text{GL}\left(2,\mathbb{C}\right)$
from $U$ and $\hat{U}$ respectively for each global section $X_{x}$.
Expand the first derivative of $f$ in a similar fashion to above,
putting
\begin{equation}
\frac{\partial f}{\partial\Omega}|=\sum_{n=-\infty}^{\infty}\phi_{n}\lambda^{n}\label{eq:expandderivative}
\end{equation}
for functions $\phi_{n}(y,z)$. We then have
\[
\hat{h}_{3}=\lambda^{-1}h_{3}\qquad\Rightarrow\qquad h_{3}=a+b\lambda
\]
and
\[
\hat{h}_{4}=h_{4}\qquad\Rightarrow\qquad h_{4}=e
\]
for functions $\left(a,b,e\right)$ on $M$. The other two components
of the splitting problem are then
\[
\hat{h}_{2}=\lambda h_{2}+\lambda\left(\sum_{n=-\infty}^{\infty}\phi_{n}\lambda^{n}\right)e
\]
\[
\Rightarrow\qquad h_{2}=-e\sum_{n=0}^{\infty}\phi_{n}\lambda^{n}
\]
and
\[
\hat{h}_{1}=h_{1}+\left(\sum_{n=-\infty}^{\infty}\phi_{n}\lambda^{n}\right)\left(a+b\lambda\right)
\]
\[
\Rightarrow\qquad h_{1}=m-\left(\sum_{n=0}^{\infty}\phi_{n}\lambda^{n}\right)\left(a+b\lambda\right)
\]
where $m$ is a new function on $M$ parametrising the non-uniqueness
in the splitting. We have $\det H=em$ so we must impose $e\neq0$
and $m\neq0$. The frame section is then
\[
v=H^{-1}\begin{pmatrix}dT|\\
d\Omega|
\end{pmatrix}=\frac{1}{em}\begin{pmatrix}e & e\sum_{n=0}^{\infty}\phi_{n}\lambda^{n}\\
-\left(a+b\lambda\right) & m-\left(\sum_{n=0}^{\infty}\phi_{n}\lambda^{n}\right)\left(a+b\lambda\right)
\end{pmatrix}\begin{pmatrix}dt-\sum_{n=1}^{\infty}d\gamma_{n}\lambda^{n}\\
dy+\lambda dz
\end{pmatrix}
\]
\[
\Rightarrow\qquad v=\frac{1}{em}\begin{pmatrix}e\left(dt+\phi_{0}dy\right)\\
m\left(dy+\lambda dz\right)-\left(a+b\lambda\right)\left(dt+\phi_{0}dy\right)
\end{pmatrix}
\]
and the clock can be read off:
\begin{equation}
\theta=m^{-1}\left(dt+\phi_{0}dy\right).\label{eq:clockin3Dtorsion}
\end{equation}
Recall that $\phi_{0}=\phi_{0}\left(y,z\right)$ and so for any choice
of $m\neq0$ we have that $d\theta\neq0$, suggesting that the moduli
space possesses Newton-Cartan torsion. The clock (\ref{eq:clockin3Dtorsion})
cannot be made compatible with any torsion-free connection as we must
have
\[
\nabla_{b}\theta_{a}=\partial_{b}\theta_{a}-\Gamma_{\,ab}^{c}\theta_{c}=0
\]
\begin{equation}
\Rightarrow\qquad\Gamma_{\,ab}^{t}=m\partial_{b}\theta_{a}-\Gamma_{\,ab}^{y}\phi_{0}.\label{eq:compatibleclock}
\end{equation}
If $\Gamma_{\,ab}^{c}$ were torsion-free then skew-symmetrising over
$ab$ in (\ref{eq:compatibleclock}) would give $d\theta=0$, which
is never true for this class of one-forms.

In stage two of the proof we construct the torsion $\Xi$-connection
of section \ref{subsec:The-Torsion--Connection} with respect to which
the clock can be made parallel. We must solve the splitting problem
\begin{equation}
\partial_{b}{\cal F}_{\,\nu}^{\mu}=-\hat{\rho}_{\,\alpha\,b}^{\mu}{\cal F}_{\,\nu}^{\alpha}+{\cal F}_{\,\beta}^{\mu}\rho_{\,\nu\,b}^{\beta}.\label{eq:this-1}
\end{equation}
in the case
\[
\partial_{b}{\cal F}_{\,\nu}^{\mu}=\partial_{b}\begin{pmatrix}1 & \frac{\partial f}{\partial\Omega}|\\
0 & \lambda^{-1}
\end{pmatrix}=\begin{pmatrix}0 & \frac{\partial^{2}f}{\partial\Omega^{2}}|\left(\delta_{b}^{y}+\delta_{b}^{z}\lambda\right)\\
0 & 0
\end{pmatrix}.
\]
Equation (\ref{eq:this-1}) constitutes four couple splitting problems;
we must find the global sections of the following four patchings:
\begin{equation}
\hat{\rho}_{\,\,T\,b}^{T}=\rho_{\,T\,b}^{T}+\frac{\partial f}{\partial\Omega}|\rho_{\,T\,b}^{\Omega}\,\,\,;\label{eq:TT}
\end{equation}
\begin{equation}
\hat{\rho}_{\,\Omega\,b}^{T}=\lambda\rho_{\,\Omega\,b}^{T}+\lambda\frac{\partial f}{\partial\Omega}|\left(\rho_{\,\Omega\,b}^{\Omega}-\hat{\rho}_{\,T\,b}^{T}\right)-\lambda\frac{\partial^{2}f}{\partial\Omega^{2}}|\left(\delta_{b}^{y}+\delta_{b}^{z}\lambda\right)\,\,\,;\label{eq:TO}
\end{equation}
\begin{equation}
\hat{\rho}_{\,T\,b}^{\Omega}=\lambda^{-1}\rho_{\,T\,b}^{\Omega}\,\,\,;\label{eq:OT}
\end{equation}
\begin{equation}
\hat{\rho}_{\,\Omega\,b}^{\Omega}=\rho_{\,\Omega\,b}^{\Omega}-\lambda\hat{\rho}_{\,T\,b}^{\Omega}\frac{\partial f}{\partial\Omega}|\,\,\,.\label{eq:OO}
\end{equation}
Equations (\ref{eq:TT}), (\ref{eq:OT}), and (\ref{eq:OO}) are immediately
tractable if we again make use of the expansion (\ref{eq:expandderivative}).
Their most general global sections are given by
\[
\rho_{\,T\,b}^{\Omega}=A_{b}^{0}+A_{b}^{1}\lambda\qquad\hat{\rho}_{\,T\,b}^{\Omega}=\hat{\lambda}A_{b}^{0}+A_{b}^{1}
\]
\[
\rho_{\,T\,b}^{T}=C_{b}-\sum_{n=0}^{\infty}\phi_{n}\lambda^{n}A_{b}^{0}-\sum_{n=-1}^{\infty}\phi_{n}\lambda^{n+1}A_{b}^{1}
\]
\[
\hat{\rho}_{\,\,T\,b}^{T}=C_{b}+\sum_{n=1}^{\infty}\phi_{-n}\hat{\lambda}^{n}A_{b}^{0}+\sum_{n=2}^{\infty}\phi_{-n}\hat{\lambda}^{n-1}A_{b}^{1}
\]
\[
\rho_{\,\Omega\,b}^{\Omega}=B_{b}+\sum_{n=1}^{\infty}\phi_{n}\lambda^{n}A_{b}^{0}+\sum_{n=0}^{\infty}\phi_{n}\lambda^{n+1}A_{b}^{1}
\]
\[
\hat{\rho}_{\,\Omega\,b}^{\Omega}=B_{b}-\sum_{n=0}^{\infty}\phi_{-n}\hat{\lambda}^{n}A_{b}^{0}-\sum_{n=1}^{\infty}\phi_{-n}\hat{\lambda}^{n-1}A_{b}^{1}
\]
where $(A^{0},A^{1},B,C)$ are arbitrary one-forms on $M$ carrying
the non-uniqueness in the splitting. The remaining equation (\ref{eq:TO})
is, after a little work, given by
\begin{multline*}
\hat{\rho}_{\,\Omega\,b}^{T}=\lambda\rho_{\,\Omega\,b}^{T}+\sum_{m=-\infty}^{\infty}\phi_{m}\lambda^{m+1}\left(B_{b}-C_{b}\right)\\
+\sum_{m=-\infty}^{\infty}\left(\left[\sum_{n=1}^{\infty}-\sum_{n=-\infty}^{-1}\right]\phi_{n}\phi_{m}\lambda^{n+m+1}A_{b}^{0}+\left[\sum_{n=0}^{\infty}-\sum_{n=-\infty}^{-2}\right]\phi_{n}\phi_{m}\lambda^{m+n+2}A_{b}^{1}\right)\\
-\lambda\frac{\partial^{2}f}{\partial\Omega^{2}}|\left(\delta_{b}^{y}+\delta_{b}^{z}\lambda\right).
\end{multline*}
This equation is the patching for an affine line bundle on $\mathbb{P}^{1}$
with underlying translation bundle ${\cal O}(-1)$ (for each direction
$x^{b}$ on $M$) and hence always has a unique solution.

If we expand the second derivative of $f$ such that
\[
\frac{\partial^{2}f}{\partial\Omega^{2}}|=\sum_{m=-\infty}^{\infty}\psi_{m}\lambda^{m}
\]
then (after a calculation) we can write the solution to the splitting
problem as
\[
\rho_{\,\Omega\,b}^{T}=-\sum_{m=0}^{\infty}\phi_{m}\lambda^{m}\left(B_{b}-C_{b}\right)-\sum_{k=1}^{\infty}\lambda^{k-1}{\cal W}_{k\,b}+\sum_{m=0}^{\infty}\left(\delta_{b}^{y}\psi_{m}+\delta_{b}^{z}\psi_{m-1}\right)\lambda^{m}
\]
\[
\hat{\rho}_{\,\Omega\,b}^{T}=\sum_{m=1}^{\infty}\phi_{-m}\hat{\lambda}^{m-1}\left(B_{b}-C_{b}\right)+\sum_{k=0}^{\infty}\hat{\lambda}^{k}{\cal W}_{-k\,b}-\sum_{m=1}^{\infty}\left(\delta_{b}^{y}\psi_{-m}+\delta_{b}^{z}\psi_{-\left(m+1\right)}\right)\hat{\lambda}^{m-1}
\]
where for convenience we define
\[
{\cal W}_{k\,b}=\sum_{n=1}^{\infty}\left[A_{b}^{0}\left(\phi_{n}\phi_{k-1-n}-\phi_{-n}\phi_{k-1+n}\right)+A_{b}^{1}\left(\phi_{n}\phi_{k-2-n}-\phi_{-n}\phi_{k-2+n}\right)\right]+A_{b}^{1}\left(\phi_{0}\phi_{k-2}+\phi_{-1}\phi_{k-1}\right).
\]

Having solved the splitting problem it is straightforward to extract
$\Gamma_{\,ab}^{c}$ from (\ref{eq:torsionconnectionreadoff}); we
have
\[
\Gamma_{\,ab}^{c}\left(\delta_{c}^{t}-\sum_{n=1}^{\infty}\left(\partial_{c}\gamma_{n}\right)\lambda^{n}\right)=-\sum_{n=1}^{\infty}\left(\partial_{a}\partial_{b}\gamma_{n}\right)\lambda^{n}+\rho_{\,T\,b}^{T}\left(\delta_{a}^{t}-\sum_{n=1}^{\infty}\left(\partial_{a}\gamma_{n}\right)\lambda^{n}\right)+\rho_{\,\Omega\,b}^{T}\left(\delta_{a}^{y}+\delta_{a}^{z}\lambda\right)
\]
and
\[
\Gamma_{\,ab}^{c}\left(\delta_{c}^{y}+\delta_{c}^{z}\lambda\right)=\rho_{\,T\,b}^{\Omega}\left(\delta_{a}^{t}-\sum_{n=1}^{\infty}\left(\partial_{a}\gamma_{n}\right)\lambda^{n}\right)+\rho_{\,\Omega\,b}^{\Omega}\left(\delta_{a}^{y}+\delta_{a}^{z}\lambda\right)
\]
from which we can read off
\[
\Gamma_{\,ab}^{y}=\delta_{a}^{t}\left[\rho_{\,T\,b}^{\Omega}\right]_{0}+\delta_{a}^{y}\left[\rho_{\,\Omega\,b}^{\Omega}\right]_{0}
\]
\[
\Gamma_{\,ab}^{z}=\delta_{a}^{t}\left[\rho_{\,T\,b}^{\Omega}\right]_{1}-\left(\partial_{a}\gamma_{1}\right)\left[\rho_{\,T\,b}^{\Omega}\right]_{0}+\delta_{a}^{y}\left[\rho_{\,\Omega\,b}^{\Omega}\right]_{1}+\delta_{a}^{z}\left[\rho_{\,\Omega\,b}^{\Omega}\right]_{0}
\]
\[
\Gamma_{\,ab}^{t}=\delta_{a}^{t}\left[\rho_{\,T\,b}^{T}\right]_{0}+\delta_{a}^{y}\left[\rho_{\,\Omega\,b}^{T}\right]_{0}
\]
where we adopt the notation $\left[\rho_{\,\nu\,b}^{\mu}\right]_{n}$
for the coefficient of $\lambda^{n}$ in $\rho_{\,\nu\,b}^{\mu}$.
The Christoffel symbols hence can be written
\[
\Gamma_{\,ab}^{y}=\delta_{a}^{t}A_{b}^{0}+\delta_{a}^{y}B_{b}
\]
\[
\Gamma_{\,ab}^{z}=\delta_{a}^{t}A_{b}^{1}+\delta_{a}^{y}\phi_{0}A_{b}^{1}+\delta_{a}^{z}\left(B_{b}-\phi_{0}A_{b}^{0}\right)
\]
\[
\Gamma_{\,ab}^{t}=\delta_{a}^{t}\left(C_{b}-\phi_{0}A_{b}^{0}-\phi_{-1}A_{b}^{1}\right)+\delta_{a}^{y}\left(-\phi_{0}\left(B_{b}-C_{b}\right)-A_{b}^{1}\phi_{0}\phi_{-1}+\delta_{b}^{y}\psi_{0}+\delta_{b}^{z}\psi_{-1}\right).
\]
This is the torsion $\Xi$-connection, a family of connections parametrised
by four arbitrary one-forms $(A^{0},A^{1},B,C)$ on $M$ which generically
possess torsion arising from the second derivative of $f$ via the
$\psi_{n}$ terms in $\Gamma_{\,ab}^{t}$. For example, 
\[
\Gamma_{\,[yz]}^{t}=\frac{1}{2}\left(-\phi_{0}\left(B_{z}-C_{z}\right)-A_{z}^{1}\phi_{0}\phi_{-1}+\psi_{-1}\right)
\]
cannot be set to zero by a global choice of $(A^{0},A^{1},B,C)$ provided
$\psi_{-1}\neq0$ and provided $\phi_{0}$ has vanishing points, which
is generically the case.

It is precisely the presence of this torsion which allows the above
connection to be made compatible with the clock (\ref{eq:clockin3Dtorsion})
in stage three of the proof. To carry out this stage we impose $\nabla\theta=0$
for the torsion $\Xi$-connection.
\[
\nabla\theta=0\qquad\Rightarrow\qquad\Gamma_{\,ab}^{t}=m\partial_{b}\theta_{a}-\Gamma_{\,ab}^{y}\phi_{0}
\]
which results in the one-form $C$ being fixed to be
\[
C_{b}=m\partial_{b}\left(m^{-1}\right)+\phi_{-1}A_{b}^{1}
\]
which simplifies $\Gamma_{\,ab}^{t}$ to 
\[
\Gamma_{\,ab}^{t}=\delta_{a}^{t}\left(m\partial_{b}\left(m^{-1}\right)-\phi_{0}A_{b}^{0}\right)+\delta_{a}^{y}\left(\phi_{0}m\partial_{b}\left(m^{-1}\right)-\phi_{0}B_{b}+\delta_{b}^{y}\psi_{0}+\delta_{b}^{z}\psi_{-1}\right).
\]
Thus the moduli space comes equipped with a family of compatible connections
and clocks with torsion parametrised by three arbitrary one-forms
$(A^{0},A^{1},B)$ and one non-vanishing function $m$.

In stage four the construction is completed with the calculation of
a family of Newton-Cartan metrics compatible with the connections
and whose kernels are spanned by the clock. The latter condition,
$h\left(\theta,\,\,\,\right)=0$, requires 
\begin{equation}
h^{at}+h^{ay}\phi_{0}=0\label{clockkernel}
\end{equation}
so it is only necessary to calculate the spatial components $h^{ij}$
of the metric; one can always reconstruct the other components by
factors of $-\phi_{0}$. Moving to the compatibility with the connection,
\begin{equation}
\nabla h=0\qquad\Rightarrow\qquad\partial_{b}h^{ij}+2h^{ij}\left(B_{b}-\phi_{0}A_{b}^{0}\right)=0.\label{eq:compatibleh}
\end{equation}
The Frobenius theorem (see, for example, \cite{MDbook}) tells us
that there exists a unique solution $h^{ij}\left(x^{b}\right)$ for
given initial data $h^{ij}\left(x_{0}^{b}\right)$ iff the one-form
$B-\phi_{0}A^{0}$ is closed. Thus we may henceforth consider $A^{0}$
to be free and $B$ to be constrained to be given by
\begin{equation}
B=\phi_{0}A^{0}+d{\cal B}\label{eq:newB}
\end{equation}
for an arbitrary function ${\cal B}$ on $M$. (We assume the first
De Rham cohomology of the relevant domain in $M$ is trivial; otherwise
(\ref{eq:newB}) is reduced to a local statement.) Equation (\ref{eq:compatibleh})
then has solutions
\[
h^{ij}=\tilde{h}^{ij}\exp\left\{ -2\int\left(B_{b}-\phi_{0}A_{b}^{0}\right)\,dx^{b}\right\} =\tilde{h}^{ij}\exp\left\{ -2{\cal B}\right\} 
\]
where $\tilde{h}^{ij}$ are constants, and for the purposes of constructing
a metric we may choose any $\tilde{h}^{ij}\in\text{GL}\left(2,\mathbb{C}\right)$.
Recall that the $h^{at}$ components can then be found from (\ref{clockkernel}).
The final form of the connection is then
\[
\Gamma_{\,ab}^{y}=\delta_{a}^{t}A_{b}^{0}+\delta_{a}^{y}\phi_{0}A_{b}^{0}+\delta_{a}^{y}\partial_{b}{\cal B}
\]
\[
\Gamma_{\,ab}^{z}=\delta_{a}^{t}A_{b}^{1}+\delta_{a}^{y}\phi_{0}A_{b}^{1}+\delta_{a}^{z}\partial_{b}{\cal B}
\]
\[
\Gamma_{\,ab}^{t}=\left(\delta_{a}^{t}+\delta_{a}^{y}\phi_{0}\right)\left(m\partial_{b}\left(m^{-1}\right)-\phi_{0}A_{b}^{0}\right)-\delta_{a}^{y}\phi_{0}\partial_{b}{\cal B}+\delta_{a}^{y}\delta_{b}^{y}\psi_{0}+\delta_{a}^{y}\delta_{b}^{z}\psi_{-1}.
\]
This completes the construction of the Newton-Cartan structures; the
free data consists of two one-forms $\left(A^{0},A^{1}\right)$, two
functions $({\cal B},m)$ (subject to $m\neq0$), and a choice of
$\tilde{h}^{ij}\in\text{GL}\left(2,\mathbb{C}\right)$.\hfill{}$\square$

Note that, compared to the case of theorem \ref{thm:conformalNCin3d},
we have a larger family of Newton-Cartan structures (depending on
more arbitrary degrees of freedom). This is because we had to use
the torsion $\Xi$-connection rather than the more powerful torsion-free
$\Lambda$-connection. It would be pleasing to be able to construct
an analogue of the $\Lambda$-connection which allows torsion; we
defer such prospects to future investigations.

\subsection{Global Vectors}

Which of the many non-relativistic conformal symmetry algebras are
singled out on $M$ by the twistor theory in the same way that Penrose's
theory singles out the conformal algebra? For the case of four dimensions
it was shown in \cite{G16} that the answer is the conformal Newton-Cartan
algebra, $\mathfrak{cnc}(3)$. In three dimensions, though, a naive
guess of $\mathfrak{cnc}(2)$ turns out to be incorrect.

In this subsection we'll identify the algebra by computing $\check{H}^{0}\left(Z,TZ\right)$
for the flat model $Z={\cal O}\oplus{\cal O}(1)$ and pulling the
elements back to $TM\rightarrow M$. The patching for $TZ\rightarrow Z$
is
\begin{equation}
\hat{\beta}^{\tilde{\alpha}}=\frac{\partial\hat{Z}^{\tilde{\alpha}}}{\partial Z^{\tilde{\beta}}}\beta^{\tilde{\beta}}\label{eq:TZpatching3d}
\end{equation}
where
\[
Z^{\tilde{\alpha}}=\left(T,\Omega,\lambda\right)^{T}\qquad\hat{Z}^{\tilde{\alpha}}=\left(\hat{T},\hat{\Omega},\hat{\lambda}\right)^{T}
\]
and where we write $\beta$ for a section of $TZ$. The components
of (\ref{eq:TZpatching3d}) read
\[
\hat{\beta}^{T}=\beta^{T}\qquad\hat{\beta}^{\Omega}=\lambda^{-1}\beta^{\Omega}-\lambda^{-2}\Omega\beta^{\lambda}\qquad\hat{\beta}^{\lambda}=-\lambda^{-2}\beta^{\lambda}.
\]
The most general global section is then given by
\[
\beta^{T}=h_{0}
\]
\[
\beta^{\Omega}=d_{0}+d_{1}\lambda+e_{0}\Omega+a_{2}\Omega\lambda+b_{1}\Omega^{2}
\]
\[
\beta^{\lambda}=a_{0}+a_{1}\lambda+a_{2}\lambda^{2}+b_{0}\Omega+b_{1}\Omega\lambda
\]
for nine arbitrary holomorphic functions $(a_{0},a_{1},a_{2},b_{0},b_{1},d_{0},d_{1},e_{0},h_{0})$
of $T$.

Now let $\tilde{X}$ be a vector field on $P\mathbb{S}^{\prime}\rightarrow M$,
so we have
\[
\tilde{X}=\tilde{X}^{\Lambda}\frac{\partial}{\partial x^{\Lambda}}=\tilde{X}^{\lambda}\frac{\partial}{\partial\lambda}+\tilde{X}^{a}\frac{\partial}{\partial x^{a}}
\]
over $U\subset\mathbb{P}^{1}$. Then $\beta=\mu_{*}\tilde{X}$ is
a vector field on $Z$ with components
\[
\beta^{\mu}=\frac{\partial Z^{\mu}|}{\partial x^{\Lambda}}\tilde{X}^{\Lambda}
\]
so we have
\[
\tilde{X}^{t}=h_{0}
\]
\[
\tilde{X}^{\lambda}=a_{0}+a_{1}\lambda+a_{2}\lambda^{2}+b_{0}\left(y+z\lambda\right)+b_{1}\left(y+z\lambda\right)\lambda
\]
\[
\Rightarrow\qquad\tilde{X}^{y}=d_{0}+e_{0}y+b_{1}y^{2}-za_{0}-b_{0}yz
\]
\[
\text{and}\qquad\tilde{X}^{z}=d_{1}+\left(e_{0}-a_{1}\right)z+a_{2}y+b_{1}yz-b_{0}z^{2}.
\]
We can then push $\tilde{X}$ down to $M$ to obtain
\[
\nu_{*}\tilde{X}=h_{0}\partial_{t}+\left(A^{i}+B_{\,j}^{i}x^{j}+x^{i}C_{j}x^{j}\right)\partial_{i}
\]
where
\[
A^{y}=d_{0}\qquad A^{z}=d_{1}
\]
\[
B_{\,j}^{i}=\begin{pmatrix}e_{0} & -a_{0}\\
a_{2} & e_{0}-a_{1}
\end{pmatrix}
\]
\[
C_{y}=b_{1}\qquad C_{z}=-b_{0}.
\]
This is a nine-dimensional Lie algebra under the usual bracket.

We can interpret this algebra heuristically as follows. Take the eight-dimensional
algebra of projective vector fields on the two-dimensional spatial
fibres, add time translations, and then promote the nine-dimensional
algebra to an infinite-dimensional one by allowing the nine components
to carry arbitrary holomorphic functions on the time axis. We can
write
\[
\check{H}^{0}\left(Z,TZ\right)=\left\{ \underset{\text{eight}}{\mathfrak{p}\left(2,\mathbb{C}\right)}\oplus\underset{\mbox{one}}{\left\{ \frac{\partial}{\partial T}\right\} }\right\} \otimes H\left({\cal O_{T}}\right)
\]
where $\mathfrak{p}\left(2,\mathbb{C}\right)$ is the algebra of projective
vector fields on the (flat) two-dimensional spatial slices, and where
$H\left({\cal O_{T}}\right)$ are the holomorphic functions on the
time axis.

\subsection{The Relativistic Limit}

Since the normal bundle $N_{x}={\cal O}\oplus{\cal O}(1)$ is stable
with respect to Kodaira deformations the Newtonian limit cannot be
realised as a jumping phenomenon in the same way as is done in \cite{DG16}
for the case of four dimensions. We can, however, carry out the reverse
procedure: if we augment the minitwistor space (with $N_{x}={\cal O}(2)$)
with an additional additive ${\cal O}(-1)$ (so that the normal bundle
is taken to be ${\cal O}(-1)\oplus{\cal O}(2)$) then one can consider
a relativistic limit of the Newtonian theory.
\begin{thm}
Let $Z_{\epsilon}\rightarrow\mathbb{P}^{1}$ be a one-parameter family
of rank-two vector bundles with patching
\[
\hat{\zeta}=\lambda\zeta+\epsilon Q
\]
\[
\hat{Q}=\lambda^{-2}Q.
\]
\end{thm}
\begin{itemize}
\item \emph{For $\epsilon=0$ $Z_{0}={\cal O}(-1)\oplus{\cal O}(2)$ and
the complex moduli space comes equipped with a three-dimensional non-degenerate
conformal structure, as in minitwistor theory.}
\item \emph{For $\epsilon\neq0$ $Z_{\epsilon}={\cal O}\oplus{\cal O}(1)$
and the complex moduli space is a three-dimensional Newton-Cartan
manifold as described in theorem \ref{thm:2+1Torsion}.}
\end{itemize}

\subsubsection*{Proof}

For $\epsilon=0$ the patching for $Z_{0}$ is exactly that of ${\cal O}(-1)\oplus{\cal O}(2)$
and so the twistor lines are
\[
\zeta|=0
\]
and
\[
Q|=\xi\lambda^{2}-2z\lambda-\tilde{\xi}
\]
which give rise to 
\[
[g]=dz^{2}+d\xi d\tilde{\xi}
\]
as expected from minitwistor theory.

For $\epsilon\neq0$ we can effect the biholomorphisms
\[
\hat{T}=\epsilon^{-1}\hat{\zeta}\qquad\hat{\Omega}=-\hat{\lambda}^{2}\epsilon^{-1}\hat{\zeta}+\hat{Q}
\]
and
\[
T=\lambda\epsilon^{-1}\zeta+Q\qquad\Omega=-\epsilon^{-1}\zeta
\]
which bring the patching to the familiar
\[
\hat{T}=T\qquad\hat{\Omega}=\lambda^{-1}\Omega
\]
and revealing that $Z_{\epsilon}={\cal O}\oplus{\cal O}(1)$ and placing
us in the arena of theorem \ref{thm:conformalNCin3d}.

\hfill{}$\square$

Note that the limit parameter $\epsilon$ appears only in the complex
structures of the one-parameter family of twistor spaces, not in the
explicit induced geometry. Thus it does not make sense to consider
the limit on the spacetime side of the correspondence in the same
way as in four dimensions \cite{DG16}.

\subsection{On Jumping Hypersurfaces of Gibbons-Hawking Manifolds\label{subsec:On-Jumping-Points}}

We end this section with a tangential result, in which three-dimensional
torsional Newton-Cartan manifolds arise on certain hypersurfaces of
Gibbons-Hawking manifolds. Recall that a Penrose twistor space is
in the Gibbons-Hawking \cite{GH} class if it admits a fibration over
${\cal O}(2)$.
\begin{thm}
Let $Z\rightarrow\mathbb{P}^{1}$ be a twistor space in the Gibbons-Hawking
class and let $\left(M,g\right)$ be its associated moduli space with
\[
g=V^{-1}\left(dt+A\right)^{2}+V\left(dz^{2}+d\xi d\tilde{\xi}\right)
\]
where the Gibbons-Hawking potential $V$ satisfies $dV=\star^{3}dA$.
On twistor lines satisfying $V=0$ the normal bundle is ${\cal O}\oplus{\cal O}(2)$
and the twistor-induced local geometry is that of a (generically-torsional)
$\left(2+1\right)$-dimensional Newton-Cartan spacetime, provided
that the restriction to $V=0$ of the flat three-metric $dz^{2}+d\xi d\tilde{\xi}$
is of rank two.
\end{thm}

\subsubsection*{Proof}

Consider first the isomorphism class of the normal bundle to twistor
lines. The patching for the normal bundle is
\[
{\cal F}=\begin{pmatrix}1 & \frac{\partial f}{\partial Q}|\\
0 & \lambda^{-2}
\end{pmatrix}
\]
and we can make an expansion
\[
\frac{\partial f}{\partial Q}|=\sum_{n=-\infty}^{\infty}\gamma_{n}\lambda^{n}.
\]
(The intersection $U\cap\hat{U}\subset\mathbb{P}^{1}$ is an annulus.)
The splitting problem is
\begin{equation}
\begin{pmatrix}1 & \frac{\partial f}{\partial Q}|\\
0 & \lambda^{-2}
\end{pmatrix}\begin{pmatrix}h_{1} & h_{2}\\
h_{3} & h_{4}
\end{pmatrix}=\begin{pmatrix}\hat{h}_{1} & \hat{h}_{2}\\
\hat{h}_{3} & \hat{h}_{4}
\end{pmatrix}\begin{pmatrix}\lambda^{m-2} & 0\\
0 & \lambda^{-m}
\end{pmatrix}\label{eq:splitting}
\end{equation}
and if there exists a holomorphic solution to this for some $m$ on
some line $X_{x}$ then the normal bundle to that line is ${\cal O}(2-m)\oplus{\cal O}(m)$.
We're interested in the exact form of $H$ on the set of twistor lines
described by $V=0$, which is when $\gamma_{-1}=0$ and where we can
solve the splitting problem for $m=2$. Henceforth assume that everything
is restricted to $\gamma_{-1}=0$.
\[
\hat{h}_{1}=h_{1}+\frac{\partial f}{\partial Q}|h_{3}
\]
\[
\hat{h}_{2}=\lambda^{2}h_{2}+\lambda^{2}\frac{\partial f}{\partial Q}|h_{4}
\]
\[
\hat{h}_{3}=\lambda^{-2}h_{3}
\]
\[
\hat{h}_{4}=h_{4}.
\]
So
\[
h_{4}=b_{0}\qquad h_{3}=a_{0}+a_{1}\lambda+a_{2}\lambda^{2}
\]
for four functions $\left(a_{0},a_{1},a_{2},b_{0}\right)$ on $M$
and
\[
h_{2}=-b_{0}\sum_{n=0}^{\infty}\gamma_{n}\lambda^{n}\qquad h_{1}=c_{0}-a_{0}\sum_{n=1}^{\infty}\gamma_{n}\lambda^{n}-a_{1}\sum_{n=0}^{\infty}\gamma_{n}\lambda^{n+1}-a_{2}\sum_{n=0}^{\infty}\gamma_{n}\lambda^{n+2}.
\]
Then
\[
\det H=\det\hat{H}=b_{0}c_{0}.
\]
One can now set
\[
a_{0}=a_{1}=a_{2}=0
\]
and
\[
b_{0}=c_{0}=1
\]
to get a simple solution to (\ref{eq:splitting}). The twistor functions
are
\[
Q|=\xi\lambda^{2}-2z\lambda-\tilde{\xi}\qquad\qquad\text{for }x^{i}=\left(\xi,\tilde{\xi},z\right)\in\mathbb{C}^{3}
\]
and
\[
T|=t-h\left(x^{i},\lambda\right)
\]
where 
\[
f|=h-\hat{h}
\]
 for $h$ and $\hat{h}$ holomorphic on $U$ and $\hat{U}$ respectively.

Now let $w=\begin{pmatrix}T\\
Q
\end{pmatrix}$ and $\hat{w}=\begin{pmatrix}\hat{T}\\
\hat{Q}
\end{pmatrix}$;
\[
d\hat{w}|={\cal F}dw|
\]
\[
\Rightarrow\qquad\hat{H}^{-1}d\hat{w}|=\begin{pmatrix}1 & 0\\
0 & \lambda^{-2}
\end{pmatrix}H^{-1}dw|
\]
is a global section of $N\otimes\Lambda^{1}\left(M\right)$ which
determines the frame $\left(\theta,e^{A^{\prime}B^{\prime}}\right)$
for the moduli space.
\[
H^{-1}dw|=\begin{pmatrix}1 & \sum_{n=0}^{\infty}\gamma_{n}\lambda^{n}\\
0 & 1
\end{pmatrix}\begin{pmatrix}dt-dh\\
d\xi\lambda^{2}-2dz\lambda-d\tilde{\xi}
\end{pmatrix}
\]
\[
\Rightarrow\qquad\begin{pmatrix}\theta\\
e^{0^{\prime}0^{\prime}}\lambda^{2}-2e^{0^{\prime}1^{\prime}}\lambda-e^{1^{\prime}1^{\prime}}
\end{pmatrix}=\begin{pmatrix}dt-dh+\left(\sum_{n=0}^{\infty}\gamma_{n}\lambda^{n}\right)\left[d\xi\lambda^{2}-2dz\lambda-d\tilde{\xi}\right]\\
d\xi\lambda^{2}-2dz\lambda-d\tilde{\xi}
\end{pmatrix}
\]
Now consider
\[
df|=dh-d\hat{h}=dQ|\sum_{n=-\infty}^{\infty}\gamma_{n}\lambda^{n}.
\]
\[
\Rightarrow\qquad dh=\left[d\xi\sum_{n=0}^{\infty}\gamma_{n}\lambda^{n+2}-2dz\sum_{n=0}^{\infty}\gamma_{n}\lambda^{n+1}-d\tilde{\xi}\sum_{n=1}^{\infty}\gamma_{n}\lambda^{n}\right]+\alpha\left[\gamma_{-2}d\xi+\gamma_{0}d\tilde{\xi}\right]
\]
for any choice of $\alpha$ (parametrising how we choose to share
this term between $dh$ and $d\hat{h}$).

The spatial part of the frame defines a conformal structure in the
two remaining spatial dimensions provided that the rank of the restriction
of $dz^{2}+d\xi d\tilde{\xi}$ to $V=0$ is of rank two, as required
in the theorem.

We can now extract the clock:
\begin{equation}
\theta=dt-\alpha\gamma_{-2}d\xi-(1-\alpha)\gamma_{0}d\tilde{\xi}\label{clock}
\end{equation}
and the triad is the standard flat triad $\left(d\xi,dz,d\hat{\xi}\right)$
restricted to $V=0$. Choose $\alpha=\frac{1}{2}$ for definiteness.
The torsion of the Newton-Cartan connection is determined by
\[
d\theta=-\frac{1}{2}d\gamma_{-2}\wedge d\xi-\frac{1}{2}d\gamma_{0}\wedge d\hat{\xi}
\]
which does not generically vanish.\hfill{}$\square$

Note that the torsion originates from $A$.

\pagebreak{}

\section{Some Novel Features in Four Dimensions\label{sec:Novel4D}}

\subsection{The $\Xi$-Connection for $Z={\cal O}\oplus{\cal O}(2)$}

The construction in this simple case amounts to taking a global section
\[
\chi_{\,\nu\,a}^{\mu}\in\check{H}^{0}\left(F|_{x},N_{x}\otimes N_{x}^{*}\otimes\Lambda_{x}^{1}\left(M\right)\right)
\]
per point $x\in M$ and extracting the connection $\Gamma_{\,bc}^{a}$
from
\[
\Gamma_{\,bc}^{a}\partial_{a}w^{\mu}|=\partial_{b}\partial_{c}w^{\mu}|+\chi_{\,\nu\,b}^{\mu}\partial_{c}w^{\nu}|+\chi_{\,\nu\,c}^{\mu}\partial_{b}w^{\nu}|.
\]
For $Z={\cal O}\oplus{\cal O}(2)$ we have
\[
N_{x}\otimes N_{x}^{*}=\begin{pmatrix}{\cal O} & {\cal O}(-2)\\
{\cal O}(2) & {\cal O}
\end{pmatrix}
\]
so the most general $\chi_{\,\nu\,a}^{\mu}$ is 
\[
\chi_{\,T\,a}^{T}=A_{a}\qquad\chi_{\,Q\,a}^{T}=0\qquad\chi_{\,Q\,a}^{Q}=E_{a}
\]
\[
\chi_{\,T\,a}^{Q}=B_{a}+\lambda C_{a}+\lambda^{2}D_{a}
\]
for five arbitrary one-forms $(A_{a},B_{a},C_{a},D_{a},E_{a})$ on
$M$. One can then read off the connection components;
\[
\Gamma_{\,tt}^{t}=2A_{t}\qquad\Gamma_{it}^{t}=A_{i}\qquad\Gamma_{\,ij}^{t}=0
\]
\[
\Gamma_{\,tt}^{\xi}=2D_{t}\qquad\Gamma_{\,tt}^{z}=-C_{t}\qquad\Gamma_{\,tt}^{\tilde{\xi}}=-2B_{t}
\]
\[
\Gamma_{\,jt}^{\xi}=D_{j}+E_{t}\delta_{j}^{\xi}\qquad\Gamma_{\,jt}^{z}=-\frac{1}{2}C_{j}+E_{t}\delta_{j}^{z}\qquad\Gamma_{\,jt}^{\tilde{\xi}}=-B_{j}+E_{t}\delta_{j}^{\tilde{\xi}}
\]
\[
\Gamma_{\,jk}^{i}=E_{j}\delta_{k}^{i}+E_{k}\delta_{j}^{i}.
\]
The connection can therefore be any connection provided that $\Gamma_{\,ij}^{t}=0$
and that the spatial sector is that of a flat projective structure
in three dimensions. Note that this includes all generalised Coriolis
forces.

\subsection{Jumps in Four-Dimensional Newtonian Twistor Theory}

Jumping phenomena were studied in \cite{DGT16} but are not constrained
to the twistor theory of (complexified) Riemannian spacetimes; they
can also occur in Newtonian twistor theory.
\begin{thm}
Let $Z$ be the total space of an affine line bundle fibred over ${\cal O}(3)$
with patching
\[
\hat{\zeta}=\lambda\zeta+f\left(S\right)
\]
\[
\hat{S}=\lambda^{-3}S
\]
where $f\left(S\right)$ is a polynomial of at least quadratic order.
$Z$ is a Newtonian twistor space: the Kodaira moduli space $M$ of
global sections is locally a complex Newton-Cartan spacetime.
\end{thm}
\begin{itemize}
\item \emph{The normal bundle is generically $N_{x}={\cal O}\oplus{\cal O}(2)$.}
\item \emph{At special points the normal bundle jumps to ${\cal O}(-1)\oplus{\cal O}(3)$;
at these jumping points the time axis in the moduli space becomes
singular.}
\end{itemize}
\emph{The special points are characterised by the vanishing of the
three two-forms induced on $M$ by the global two form $d\hat{\zeta}\wedge d\hat{S}$,
or equivalently by the vanishing of the constant term $\gamma_{0}$
in the Laurent expansion of $\frac{\partial f}{\partial S}|$.}

\subsubsection*{Proof}

Write the global sections of ${\cal O}(3)$ as
\[
S|=x_{0}+x_{1}\lambda+x_{2}\lambda^{2}+x_{3}\lambda^{3}.
\]
First we will consider the normal bundle; we will show that generically
the isomorphism class of $N_{x}$ is ${\cal O}\oplus{\cal O}(2)$.
The patching for $N_{x}$ is
\[
{\cal F}=\begin{pmatrix}\lambda & \frac{\partial f}{\partial S}|\\
0 & \lambda^{-3}
\end{pmatrix}
\]
and the splitting problem to be solved is 
\begin{equation}
\hat{h}_{1}=\lambda h_{1}+\frac{\partial f}{\partial S}|h_{3}\label{hh1-4}
\end{equation}
\begin{equation}
\hat{h}_{2}=\lambda^{3}h_{2}+\lambda^{2}\frac{\partial f}{\partial S}|h_{4}\label{hh2-4}
\end{equation}
\[
\hat{h}_{3}=\lambda^{-3}h_{3}
\]
\[
\hat{h}_{4}=\lambda^{-1}h_{4}.
\]
As usual we put
\[
h_{3}=\sum_{n=0}^{3}a_{n}\lambda^{n}\qquad\text{and}\qquad h_{4}=b_{0}+b_{1}\lambda.
\]
To proceed further we make an expansion:
\[
\frac{\partial f}{\partial S}|=\gamma_{0}+\sum_{n=1}^{\infty}\gamma_{n}\lambda^{n}
\]
where $\gamma_{0}$ is a function of $x_{0}$ only. For a global solution
of (\ref{hh2-4}) we then require
\begin{equation}
b_{0}\gamma_{0}=0\,\,,\label{eq:jNTTconstraint}
\end{equation}
leaving 
\[
\hat{h}_{2}=0.
\]
We also have
\[
\hat{h}_{1}=\gamma_{0}a_{0}.
\]
The determinant is then
\[
\det\hat{H}=\gamma_{0}a_{0}b_{1}\,\,,
\]
and we conclude that $N_{x}={\cal O}\oplus{\cal O}(2)$ for $\gamma_{0}\neq0$.
A straightforward calculation then shows that $N_{x}={\cal O}(-1)\oplus{\cal O}(3)$
for any twistor line $X_{x}$ with $\gamma_{0}=0$.

Since we generically have $N_{x}={\cal O}\oplus{\cal O}(2)$ we anticipate
that the moduli space will be a Newton-Cartan spacetime; to see this
in detail we must construct the twistor functions. In this case it
is very straightforward: over $\hat{U}$ we have
\[
\hat{\zeta}=f\left(S|\left(\lambda=0\right)\right):={\cal T}(x_{0})
\]
and
\[
\hat{S}|=x_{0}\hat{\lambda}^{3}+x_{1}\hat{\lambda}^{2}+x_{2}\hat{\lambda}+x_{3}.
\]
We note that $\gamma_{0}=\frac{d{\cal T}}{dx_{0}}$. The geometry
induced on the moduli space can be found by identifying null vectors
as those tangent to alpha surfaces. The clock is therefore
\[
\theta=\alpha\left(x_{0}\right)\gamma_{0}\left(x_{0}\right)dx_{0}
\]
for any non-vanishing $\alpha$ and the conformal covariant Galilean
metric is 
\[
h^{-1}=\beta\left(dx_{2}^{2}-4dx_{1}dx_{3}\right)
\]
for any non-vanishing $\beta$. The geometry is therefore Newton-Cartan
at generic points. At the jumping points with $\gamma_{0}=0$, though,
the clock vanishes.

One can construct a map
\begin{equation}
x_{0}\mapsto t(x_{0})={\cal T}(x_{0})\label{maptostandardNC}
\end{equation}
taking $M$ to a more usual non-jumping Newton-Cartan spacetime but
the map is not a diffeomorphism; it eliminates the jumping points.
This situation is entirely analogous to the map taking the jumping
spacetimes from \cite{DGT16} to Gibbons-Hawking form.

The three two-forms arising from the restriction to twistor lines
of
\[
d\hat{\zeta}\wedge d\hat{S}=\lambda^{-2}d\zeta\wedge dS
\]
clearly vanish on and only on the jumping points.\hfill{}$\square$

\subsubsection*{Example}
\begin{flushleft}
The simplest example of a jumping Newtonian twistor space is when
$f\left(S\right)=\frac{1}{2}S^{2}$. The conformal clock admits a
representative
\[
\theta=x_{0}dx_{0}\,\,,
\]
vanishing at one point $x_{0}=0$. The map (\ref{maptostandardNC})
is
\[
t=\frac{1}{2}x_{0}^{2}\,\,\,,
\]
and so the Newton-Cartan spacetime is thus a $2$-fold cover of the
standard Newton-Cartan spacetime (with time coordinate $t$), branched
over the spatial fibre $t=0$.
\par\end{flushleft}

\pagebreak{}

\section{Five Dimensions\label{sec:Five-Dimensions}}

In this section we'll study some five-dimensional Kodaira families.
Mostly we'll be concerned with the Newtonian theory, for which the
normal bundle will be ${\cal O}\oplus{\cal O}(1)\oplus{\cal O}(1)$,
though later in this section we will also consider the normal bundle
${\cal O}(4)$, for which the moduli space comes equipped with relativistic
geometry in accordance with theorem \ref{thm:O(2n)}.

\subsection{Galilean Structures and Canonical Connections}

We'll begin by constructing the geometry on the moduli space for the
undeformed Newtonian case, equipped with its canonical $\Lambda$-connection.
\begin{thm}
Let $Z={\cal O}\oplus{\cal O}(1)\oplus{\cal O}(1)$. The induced geometry
on the moduli space of global sections of $Z\rightarrow\mathbb{P}^{1}$
is that of a complex five-dimensional Newton-Cartan spacetime $\left(M,h,\theta,\nabla\right)$,
where the connection components of $\nabla$ depend (up to diffeomorphisms
of the time axis) on one conformal factor and seven arbitrary functions,\label{thm:flat5d}
\end{thm}
\begin{itemize}
\item \emph{four of which are the Newtonian gravitational force $\Gamma_{\,tt}^{i}$;}
\item \emph{and the remaining three of form an anti-self-dual spatial two-form
$W_{ij}$ describing Coriolis forces.}
\end{itemize}

\subsubsection*{Proof}

The normal bundle to twistor lines is $N_{x}={\cal O}\oplus{\cal O}(1)\oplus{\cal O}(1)$,
which satisfies $\check{H}^{1}\left(\mathbb{P}^{1},N_{x}\right)=0$.
By the Kodaira theorem \cite{Kodaira} the moduli space $M$ is therefore
a complex manifold, with dimension $\text{dim}\check{H}^{0}\left(\mathbb{P}^{1},N_{x}\right)=5$.

To construct the conformal Galilean structure $\left(h,\theta\right)$
we first find the twistor lines explicitly. In homogeneous coordinates
the patching is
\[
\hat{T}=T\qquad\mbox{(weight zero)}
\]
and
\[
\hat{\omega}^{A}=\omega^{A}\qquad\mbox{(weight one)},
\]
with twistor lines
\[
T|=t\qquad\text{and}\qquad\omega^{A}|=x^{AA^{\prime}}\pi_{A^{\prime}}
\]
for coordinates $\left(t,x^{AA^{\prime}}\right)$ on $M$. To find
the frame we need to solve
\[
\begin{pmatrix}1 & 0 & 0\\
0 & \lambda^{-1} & 0\\
0 & 0 & \lambda^{-1}
\end{pmatrix}\begin{pmatrix}h_{1} & h_{2} & h_{3}\\
h_{4} & h_{5} & h_{6}\\
h_{7} & h_{8} & h_{9}
\end{pmatrix}=\begin{pmatrix}\hat{h}_{1} & \hat{h}_{2} & \hat{h}_{3}\\
\hat{h}_{4} & \hat{h}_{5} & \hat{h}_{6}\\
\hat{h}_{7} & \hat{h}_{8} & \hat{h}_{9}
\end{pmatrix}\begin{pmatrix}1 & 0 & 0\\
0 & \lambda^{-1} & 0\\
0 & 0 & \lambda^{-1}
\end{pmatrix}.
\]
for the most general $H$ and $\hat{H}$. The solution is
\[
H^{-1}=\begin{pmatrix}\frac{1}{m} & 0 & 0\\
\frac{1}{mk}\sum_{p=0}^{1}\left(k_{2}b_{p}-k_{4}a_{p}\right)\lambda^{p} & \frac{k_{4}}{k} & -\frac{k_{2}}{k}\\
\frac{1}{mk}\sum_{p=0}^{1}\left(k_{3}a_{p}-k_{1}b_{p}\right)\lambda^{p} & -\frac{k_{3}}{k} & \frac{k_{1}}{k}
\end{pmatrix}.
\]
for nine arbitrary functions $(m,a_{0},a_{1},b_{0},b_{1},k_{1},k_{2},k_{3},k_{4})$
on $M$ which must be chosen such that $m\neq0$ and $k:=\left(k_{1}k_{4}-k_{2}k_{3}\right)\neq0$
anywhere. These functions will determine the ambiguity in the moduli
space geometry. The frame section is therefore given by
\[
v=\begin{pmatrix}\frac{1}{m} & 0 & 0\\
\frac{1}{mk}\sum_{p=0}^{1}\left(k_{2}b_{p}-k_{4}a_{p}\right)\lambda^{p} & \frac{k_{4}}{k} & -\frac{k_{2}}{k}\\
\frac{1}{mk}\sum_{p=0}^{1}\left(k_{3}a_{p}-k_{1}b_{p}\right)\lambda^{p} & -\frac{k_{3}}{k} & \frac{k_{1}}{k}
\end{pmatrix}\begin{pmatrix}dt\\
dx^{01^{\prime}}+dx^{00^{\prime}}\lambda\\
dx^{11^{\prime}}+dx^{10^{\prime}}\lambda
\end{pmatrix}
\]
from which we can read-off the frame
\[
\theta=m^{-1}\,dt
\]
\[
e^{00^{\prime}}=\frac{k_{4}}{k}dx^{00^{\prime}}-\frac{k_{2}}{k}dx^{10^{\prime}}+\frac{1}{mk}\left(k_{2}b_{1}-k_{4}a_{1}\right)dt
\]
\[
e^{01^{\prime}}=\frac{k_{4}}{k}dx^{01^{\prime}}-\frac{k_{2}}{k}dx^{11^{\prime}}+\frac{1}{mk}\left(k_{2}b_{0}-k_{4}a_{0}\right)dt
\]
\[
e^{10^{\prime}}=-\frac{k_{3}}{k}dx^{00^{\prime}}+\frac{k_{1}}{k}dx^{10^{\prime}}+\frac{1}{mk}\left(k_{3}a_{1}-k_{1}b_{1}\right)dt
\]
\[
e^{11^{\prime}}=-\frac{k_{3}}{k}dx^{01^{\prime}}+\frac{k_{1}}{k}dx^{11^{\prime}}+\frac{1}{mk}\left(k_{3}a_{0}-k_{1}b_{0}\right)dt\,,
\]
leading to the natural decomposition of the tangent bundle
\[
TM=\mathbb{C}\oplus\left(\mathbb{S}\otimes\mathbb{S}^{\prime}\right).
\]
The conformal clock is therefore $\theta=m^{-1}\,dt$, and the conformal
\emph{covariant} metric is
\[
h^{-1}=\epsilon_{AB}\epsilon_{A^{\prime}B^{\prime}}e_{\,}^{AA^{\prime}}\otimes e_{\,}^{BB^{\prime}}
\]
\begin{multline*}
\Rightarrow\qquad h^{-1}=k^{-1}\left(dx^{00^{\prime}}dx^{11^{\prime}}-dx^{10^{\prime}}dx^{01^{\prime}}\right)\\
+m^{-1}k^{-1}dt\left(-b_{0}dx^{00^{\prime}}+b_{1}dx^{01^{\prime}}+a_{0}dx^{10^{\prime}}-a_{1}dx^{11^{\prime}}\right)\\
+m^{-2}k^{-1}\left(a_{1}b_{0}-a_{0}b_{1}\right)dt^{2}
\end{multline*}
The metric is rank-four, as one would expect for a five-dimensional
Newton-Cartan spacetime. The contravariant metric is obtained as the
projective inverse in the following way. First find a vector $U$
such that
\[
\theta(U)=1\qquad\text{and}\qquad h^{-1}(U,\,\,)=0\,\,,
\]
which uniquely determines
\[
U=m\partial_{t}+b_{0}\partial_{11^{\prime}}+b_{1}\partial_{10^{\prime}}+a_{0}\partial_{01^{\prime}}+a_{1}\partial_{00^{\prime}}.
\]
The contravariant metric $h$ is then the unique solution to
\[
h^{ab}h_{bc}=\delta_{c}^{a}-U^{a}\theta_{c}\qquad\text{and}\qquad h(\theta,\,\,)=0\,,
\]
which determines
\[
h=k^{-1}\epsilon^{AB}\epsilon^{A^{\prime}B^{\prime}}\frac{\partial}{\partial x^{AA^{\prime}}}\otimes\frac{\partial}{\partial x^{BB^{\prime}}}.
\]
Thus we have a Galilean structure $\left(M,h,\theta\right)$ depending
on arbitrary functions. (We could also equivalently have used the
more traditional twistor theory method of calculating the null vectors
from the twistor functions as was done for the case of four dimensions
in \cite{DG16}.) 

It only remains to calculate the physical induced connection, i.e.
we must construct the $\Lambda$-connection. Denote by 
\begin{equation}
\hat{w}^{\mu}=\begin{pmatrix}\hat{T}\\
\nicefrac{\hat{\omega}^{A}}{\pi_{0^{\prime}}}
\end{pmatrix}\quad\mbox{and}\quad w^{\mu}=\begin{pmatrix}T\\
\nicefrac{\omega^{A}}{\pi_{1^{\prime}}}
\end{pmatrix}\label{eq:inhomNTT5coords}
\end{equation}
column vectors of inhomogeneous twistor coordinates on the fibres,
and for ease of notation set
\begin{equation}
x^{AA^{\prime}}=x^{i}=\begin{pmatrix}v & u\\
y & x
\end{pmatrix}.\label{eq:spacecoordfor5d}
\end{equation}

Following the discussion in section \ref{subsec:Induced-Affine-Connections}
the construction of $\nabla$ occurs in two stages, the first being
the solution of the splitting problem
\begin{equation}
{\cal F}_{\,\nu\rho}^{\mu}=-\hat{\sigma}_{\,\alpha\beta}^{\mu}{\cal F}_{\,\nu}^{\alpha}{\cal F}_{\,\rho}^{\beta}+{\cal F}_{\,\gamma}^{\mu}\sigma_{\,\nu\rho}^{\gamma}\label{eq:splitting-1}
\end{equation}
for a $0$-cochain $\{\sigma\}$ of $N\otimes\left(N^{*}\odot N^{*}\right)\rightarrow\mathbb{P}^{1}$,
where
\[
{\cal F}_{\,\nu}^{\alpha}=\frac{\partial\hat{w}^{\alpha}}{\partial w^{\nu}}|\quad\mbox{and}\quad{\cal F}_{\,\nu\rho}^{\mu}=\frac{\partial^{2}\hat{w}^{\mu}}{\partial w^{\nu}\partial w^{\rho}}|.
\]
The solution of (\ref{eq:splitting-1}) depends on nine arbitrary
functions (on $M$) because
\[
\check{H}^{0}\left(\mathbb{P}^{1},N\otimes\left(N^{*}\odot N^{*}\right)\right)=\mathbb{C}^{9}\,\,,
\]
and is explicitly given by
\[
\hat{\sigma}_{\,AB}^{T}=\sigma_{\,AB}^{T}=0\qquad\hat{\sigma}_{\,AT}^{T}=\sigma_{\,AT}^{T}=0\qquad\hat{\sigma}_{\,BC}^{A}=\sigma_{\,BC}^{A}=0
\]
\[
\hat{\sigma}_{\,TT}^{T}=\sigma_{\,TT}^{T}=\Sigma\qquad\hat{\sigma}_{\,BT}^{A}=\sigma_{\,BT}^{A}=\chi_{\,B}^{A}
\]
\[
\hat{\sigma}_{\,TT}^{A}=\lambda^{-1}\phi^{A}+\psi^{A}\qquad\sigma_{\,TT}^{A}=\phi^{A}+\lambda\psi^{A}.
\]
(The nine functions are $\Sigma$, $\phi^{A}$, $\psi^{A}$, and $\chi_{\,B}^{A}$.)

The second stage of the construction is the reading-off of $\Gamma_{\,bc}^{a}$
from the map $T^{[2]}M\rightarrow TM$ determined by $\{\sigma\}$
via the Kodaira isomorphism $TM=\check{H}^{0}\left(\mathbb{P}^{1},N\right)$.
Concretely we read off $\Gamma_{\,bc}^{a}(x^{d})$ from 
\[
\Gamma_{\,bc}^{a}\partial_{a}w^{\mu}|=\partial_{b}\partial_{c}w^{\mu}|+\sigma_{\,\nu\rho}^{\mu}\partial_{b}w^{\nu}|\partial_{c}w^{\rho}|\,\,,
\]
giving us
\begin{equation}
\Gamma_{\,tt}^{t}=\Sigma\label{merk1}
\end{equation}
\[
\Gamma_{\,tt}^{u}=\phi^{0}\quad\Gamma_{\,tt}^{v}=\psi^{0}\quad\Gamma_{\,tt}^{x}=\phi^{1}\quad\Gamma_{\,tt}^{y}=\psi^{1}
\]
\[
\Gamma_{\,ut}^{u}=\Gamma_{\,vt}^{v}=\chi_{\,0}^{0}\quad\Gamma_{\,xt}^{x}=\Gamma_{\,yt}^{y}=\chi_{\,1}^{1}
\]
\begin{equation}
\Gamma_{\,xt}^{u}=\Gamma_{\,yt}^{v}=\chi_{\,1}^{0}\quad\Gamma_{\,ut}^{x}=\Gamma_{\,vt}^{y}=\chi_{\,0}^{1}\label{merk4}
\end{equation}
with all other components of $\Gamma_{\,bc}^{a}$ vanishing. Two of
these functions are related to ones we already have by the compatibility
conditions
\[
\nabla\theta=0\qquad\nabla h=0\,\,,
\]
which give us
\[
\Sigma=-\partial_{t}\ln m\quad\mbox{and}\quad\mbox{tr}(\chi)=\chi_{\,0}^{0}+\chi_{\,1}^{1}=-\frac{1}{2}\partial_{t}\ln k\,\,,
\]
as well as
\[
\frac{\partial m}{\partial x^{AA^{\prime}}}=\frac{\partial k}{\partial x^{AA^{\prime}}}=0
\]
so these two factors are functions of time only. The function $m$
can be set to one without loss of generality by a diffeomorphism of
the time axis.

The four components $\Gamma_{\,tt}^{i}$ are completely arbitrary
(given in terms of $\phi^{A}$ and $\psi^{A})$, whilst the remaining
$\Gamma_{\,jt}^{i}$ components depend on three arbitrary functions
from the traceless part of $\chi_{\,B}^{A}$. Thus we get only three
functions' worth of $\Gamma_{\,jt}^{i}$ instead of the most general
case depending on six functions. Given that this is twistor theory
it is perhaps no surprise that the three functions form an anti-self-dual
two-form on spatial fibres. Concretely we have 
\[
\Gamma_{\,jt}^{i}=\delta^{ik}W_{jk}
\]
 with
\begin{equation}
W=\left(\chi_{\,0}^{0}-\chi_{\,1}^{1}\right)\left[du\wedge dy+dx\wedge dv\right]+2\chi_{\,1}^{0}\left[dx\wedge dy\right]+2\chi_{\,0}^{1}\left[dv\wedge du\right].\label{asd2form}
\end{equation}
It is then straightforward to check that (\ref{asd2form}) is the
most general anti-self-dual two-form on spatial fibres with respect
to a volume-form
\[
\epsilon_{\mbox{space}}=dx\wedge dy\wedge dv\wedge du\,\,,
\]
completing the proof.

\begin{flushright}
$\square$
\par\end{flushright}

We thus conclude that Newtonian twistor theory in five dimensions
admits generalised Coriolis connection components as arbitrary functions
in its $\Lambda$-connection, but only \emph{half} of them.

\subsubsection*{$\Xi$-Connection}

As described in section \ref{subsec:Induced-Affine-Connections} the
calculation amounts to taking a global section $\chi_{\,\nu\,a}^{\mu}$
of $N_{x}\otimes N_{x}^{*}\otimes\Lambda_{x}^{1}\left(M\right)$ per
point $x\in M$. For $N_{x}={\cal O}\oplus{\cal O}(1)\oplus{\cal O}(1)$
we have 
\[
N_{x}\otimes N_{x}^{*}=\begin{pmatrix}{\cal O} & {\cal O}(-1) & {\cal O}(-1)\\
{\cal O}(1) & {\cal O} & {\cal O}\\
{\cal O}(1) & {\cal O} & {\cal O}
\end{pmatrix}
\]
giving us
\[
\chi_{\,T\,a}^{T}=C_{a}\qquad\chi_{\,A\,a}^{T}=0\qquad\chi_{\,T\,a}^{A}=A_{\,a}^{AA^{\prime}}\pi_{A^{\prime}}
\]
\[
\chi_{\,B\,a}^{A}=B_{\,B\,a}^{A}
\]
for nine arbitrary one-forms $\left(C,A^{AA^{\prime}},B_{\,B}^{A}\right)$
on $M$ constituting forty-five arbitrary functions. We extract the
connection symbols by reading off from
\[
\Gamma_{\,bc}^{a}\partial_{a}w^{\mu}|=\partial_{b}\partial_{c}w^{\mu}|+\chi_{\,\nu\,b}^{\mu}\partial_{c}w^{\nu}|+\chi_{\,\nu\,c}^{\mu}\partial_{b}w^{\nu}|\,\,,
\]
which gives us 
\[
\Gamma_{\,bc}^{t}=\delta_{\,b}^{t}C_{c}+\delta_{\,c}^{t}C_{b}
\]
\[
\Gamma_{\,\,\,tt}^{AA^{\prime}}\pi_{A^{\prime}}=2A_{\,t}^{AB^{\prime}}\pi_{B^{\prime}}
\]
\[
\Gamma_{\,\,\,BB^{\prime}t}^{AA^{\prime}}\pi_{A^{\prime}}=A_{\,BB^{\prime}}^{AC^{\prime}}\pi_{C^{\prime}}+B_{\,B\,t}^{A}\pi_{B^{\prime}}
\]
\[
\Gamma_{\,\,\,BB^{\prime}CC^{\prime}}^{AA^{\prime}}\pi_{A^{\prime}}=B_{\,C\,BB^{\prime}}^{A}\pi_{C^{\prime}}+B_{\,B\,CC^{\prime}}^{A}\pi_{B^{\prime}}.
\]
The $\Xi$-connection therefore contains every connection which has
$\Gamma_{\,ij}^{t}=0$ and for which the four-dimensional spatial
sector resembles that of the $\Xi$-connection (\ref{eq:XiforO1O1})
for ${\cal O}(1)\oplus{\cal O}(1)$ as discussed in \ref{subsec:Induced-Affine-Connections}.

\subsection{Deformations in the Newtonian Theory\label{subsec:FiveDDeformations}}

In this subsection we'll study deformations of the form
\begin{equation}
\hat{T}=T+\epsilon f\left(\Omega^{A},\lambda\right)\label{eq:def5}
\end{equation}
over the total space of ${\cal O}(1)\oplus{\cal O}(1)$, where $\epsilon$
is a deformation parameter. Of course, from one point of view this
is nothing more than a Ward bundle on the twistor space for flat four-dimensional
spacetime. The approach adopted in this paper is instead to study
the geometry of the full five-dimensional moduli space of global sections.

Given that the normal bundle $N_{x}={\cal O}\oplus{\cal O}(1)\oplus{\cal O}(1)$
is stable with respect to all Kodaira deformations we must therefore
possess a five-dimensional Galilean structure, but the connection
is more subtle. $\check{H}^{1}\left(\mathbb{P}^{1},N_{x}\otimes\left(N_{x}^{*}\odot N_{x}^{*}\right)\right)\neq0$,
so some deformations will result in a moduli space which does \emph{not}
possess a $\Lambda$-connection. Like in the three-dimensional case,
what is going wrong is that a $\Lambda$-connection is, by construction,
torsion-free, and deformations of the form (\ref{eq:def5}) give rise
to moduli spaces whose Newton-Cartan structures possess torsion.
\begin{thm}
\label{thm:5Dtorsion}Let $Z$ be a complex four-fold fibred over
$\mathbb{P}^{1}$ with patching given by (\ref{eq:def5}) whose five-parameter
family of global sections $X_{x}$ have normal bundle $X_{x}={\cal O}\oplus{\cal O}(1)\oplus{\cal O}(1)$.
The moduli space $M$ of those sections is a complex five-dimensional
manifold equipped with a Galilean structure with torsion whose clock
admits a representative with
\[
d\theta=\epsilon\left\{ \frac{1}{2\pi i}\oint\frac{\partial^{2}f}{\partial\omega^{A}\partial\omega^{B}}|\frac{\pi_{B^{\prime}}}{\pi_{0^{\prime}}}\pi\cdot d\pi\right\} dx^{BB^{\prime}}\wedge dx^{A1^{\prime}}.
\]
\end{thm}
(Recall that $\omega^{A}$ are the homogeneous versions of $\Omega^{A}$.)

\subsubsection*{Proof}

Take the global sections of the base ${\cal O}(1)\oplus{\cal O}(1)$
to be $x^{AA^{\prime}}\pi_{A^{\prime}}$ as in (\ref{eq:spacecoordfor5d}),
or in inhomogeneous coordinates
\[
\Omega^{0}|=u+v\lambda\qquad\Omega^{1}|=x+y\lambda.
\]
Now restrict $f$ to these lines and expand it in a Laurent series
in $\lambda$;
\[
f|=\sum_{n=-\infty}^{\infty}\gamma_{n}\lambda^{n}\qquad\text{for}\qquad\gamma_{n}=\frac{1}{2\pi i}\oint f\left(\Omega^{A}|,\lambda\right)\lambda^{-\left(1+n\right)}d\lambda.
\]
The global sections of $Z\rightarrow\mathbb{P}^{1}$ are then completed
by
\[
T|=t-\epsilon\sum_{n=1}^{\infty}\gamma_{n}\lambda^{n}.
\]
The next task is to calculate the frame section. We must solve
\[
\begin{pmatrix}1 & \epsilon\frac{\partial f}{\partial\Omega^{0}}| & \epsilon\frac{\partial f}{\partial\Omega^{1}}|\\
0 & \lambda^{-1} & 0\\
0 & 0 & \lambda^{-1}
\end{pmatrix}\begin{pmatrix}h_{1} & h_{2} & h_{3}\\
h_{4} & h_{5} & h_{6}\\
h_{7} & h_{8} & h_{9}
\end{pmatrix}=\begin{pmatrix}\hat{h}_{1} & \hat{h}_{2} & \hat{h}_{3}\\
\hat{h}_{4} & \hat{h}_{5} & \hat{h}_{6}\\
\hat{h}_{7} & \hat{h}_{8} & \hat{h}_{9}
\end{pmatrix}\begin{pmatrix}1 & 0 & 0\\
0 & \lambda^{-1} & 0\\
0 & 0 & \lambda^{-1}
\end{pmatrix}.
\]
Written out in full we must therefore solve
\[
\hat{h}_{1}=h_{1}+\epsilon\frac{\partial f}{\partial\Omega^{0}}|h_{4}+\epsilon\frac{\partial f}{\partial\Omega^{1}}|h_{7}
\]
\[
\hat{h}_{2}=\lambda h_{2}+\epsilon\lambda\frac{\partial f}{\partial\Omega^{0}}|h_{5}+\epsilon\lambda\frac{\partial f}{\partial\Omega^{1}}|h_{8}
\]
\[
\hat{h}_{3}=\lambda h_{3}+\epsilon\lambda\frac{\partial f}{\partial\Omega^{0}}|h_{6}+\epsilon\lambda\frac{\partial f}{\partial\Omega^{1}}|h_{9}
\]
\[
\hat{h}_{4}=\lambda^{-1}h_{4}\qquad\hat{h}_{7}=\lambda^{-1}h_{7}
\]
\[
\hat{h}_{5}=h_{5}\qquad\hat{h}_{6}=h_{6}\qquad\hat{h}_{8}=h_{8}\qquad\hat{h}_{9}=h_{9}.
\]
Put
\[
\begin{pmatrix}h_{5} & h_{6}\\
h_{8} & h_{9}
\end{pmatrix}=\begin{pmatrix}k_{1} & k_{2}\\
k_{3} & k_{4}
\end{pmatrix}
\]
and
\[
\begin{pmatrix}h_{4}\\
h_{7}
\end{pmatrix}=\begin{pmatrix}a_{0}+a_{1}\lambda\\
b_{0}+b_{1}\lambda
\end{pmatrix}.
\]
Now expand the derivatives of $f$ (restricted to twistor lines) in
Laurent series;
\[
\frac{\partial f}{\partial\Omega^{A}}|=\sum_{n=-\infty}^{\infty}\phi_{n,A}\lambda^{n}\qquad\text{for}\qquad\phi_{n,A}=\frac{1}{2\pi i}\oint\frac{\partial f}{\partial\Omega^{A}}|\lambda^{-\left(1+n\right)}d\lambda.
\]
We then (uniquely) obtain
\[
h_{2}=-\epsilon\sum_{n=0}^{\infty}\left(\phi_{n,0}k_{1}+\phi_{n,1}k_{3}\right)\lambda^{n}
\]
\[
h_{3}=-\epsilon\sum_{n=0}^{\infty}\left(\phi_{n,0}k_{2}+\phi_{n,1}k_{4}\right)\lambda^{n}
\]
and we can solve for the remaining piece $h_{1}$ up to the arbitrary
function $m$ to obtain
\[
h_{1}=m-\epsilon\sum_{n=0}^{\infty}\lambda^{n}\left[\phi_{n,0}\left(a_{0}+a_{1}\lambda\right)+\phi_{n,1}\left(b_{0}+b_{1}\lambda\right)\right].
\]
Define $k=k_{1}k_{4}-k_{2}k_{3}$. The determinant of $H$ is then
given by
\[
\det H=mk
\]
so we must impose $m\neq0$ and $k\neq0$. 

We can then calculate
\[
\left(H^{-1}\right)_{\,T}^{T}=m^{-1}
\]
\[
\left(H^{-1}\right)_{\,A}^{T}=\epsilon m^{-1}\sum_{n=0}^{\infty}\phi_{n,A}\lambda^{n}
\]
\[
\left(H^{-1}\right)_{\,T}^{A}=\begin{pmatrix}\frac{1}{mk}\sum_{p=0}^{1}\left(k_{2}b_{p}-k_{4}a_{p}\right)\lambda^{p}\\
\frac{1}{mk}\sum_{p=0}^{1}\left(k_{3}a_{p}-k_{1}b_{p}\right)\lambda^{p}
\end{pmatrix}
\]
\[
\left(H^{-1}\right)_{\,B}^{A}=\begin{pmatrix}\frac{k_{4}}{k}+\frac{\epsilon}{mk}\sum_{p=0}^{1}\sum_{n=0}^{\infty}\phi_{n,0}\left(k_{2}b_{p}-k_{4}a_{p}\right)\lambda^{n+p} & -\frac{k_{2}}{k}+\frac{\epsilon}{mk}\sum_{p=0}^{1}\sum_{n=0}^{\infty}\phi_{n,1}\left(k_{2}b_{p}-k_{4}a_{p}\right)\lambda^{n+p}\\
-\frac{k_{3}}{k}+\frac{\epsilon}{mk}\sum_{p=0}^{1}\sum_{n=0}^{\infty}\phi_{n,0}\left(k_{3}a_{p}-k_{1}b_{p}\right)\lambda^{n+p} & \frac{k_{1}}{k}+\frac{\epsilon}{mk}\sum_{p=0}^{1}\sum_{n=0}^{\infty}\phi_{n,1}\left(k_{3}a_{p}-k_{1}b_{p}\right)\lambda^{n+p}
\end{pmatrix}.
\]
The clock is therefore given by
\[
\theta=m^{-1}dT|+m^{-1}\epsilon\sum_{n=0}^{\infty}\phi_{n,A}\lambda^{n}d\Omega^{A}|
\]
\[
\Rightarrow\qquad\theta=m^{-1}\left(dt-\epsilon\sum_{n=1}^{\infty}d\gamma_{n}\lambda^{n}+\epsilon\sum_{n=0}^{\infty}\phi_{n,A}\lambda^{n}d\Omega^{A}|\right).
\]
Now, we have
\[
d\gamma_{n}=\phi_{n,A}dx^{A1^{\prime}}+\phi_{n-1,A}dx^{A,0^{\prime}}
\]
so
\[
\theta=m^{-1}\left(dt+\epsilon\phi_{0,A}dx^{A1^{\prime}}\right).
\]
As in the five-dimensional case this is not closed for any $m\neq0$
(and $\epsilon\neq0$). Now take a representative with $m=1$; we
then have
\[
d\theta=\epsilon\partial_{BB^{\prime}}\phi_{0,A}dx^{BB^{\prime}}\wedge dx^{A1^{\prime}}
\]
\[
\Rightarrow\qquad d\theta=\epsilon\left\{ \frac{1}{2\pi i}\oint\frac{\partial^{2}f}{\partial\omega^{A}\partial\omega^{B}}|\frac{\pi_{B^{\prime}}}{\pi_{0^{\prime}}}\pi\cdot d\pi\right\} dx^{BB^{\prime}}\wedge dx^{A1^{\prime}}\,\,.
\]
The Newton-Cartan metric arises, as in theorem \ref{thm:flat5d},
from the projective inverse of degenerate covariant metric arising
from the frame section, completing the construction of a Galilean
structure with torsion.

\hfill{}$\square$

In theorem \ref{thm:5Dtorsion} we chose to merely construct the torsional
Galilean structure, exhibiting the torsion via the non-closure of
the clock. We could, however, go further and explicitly construct
the torsion $\Xi$-connection of section \ref{subsec:The-Torsion--Connection}
as was done for the three-dimensional case in theorem \ref{thm:2+1Torsion}.
In the interests of brevity we omit this cumbersome calculation.

\subsection{Global Vectors}

In this section we will study the image on $M$ of $\check{H}^{0}\left(Z,TZ\right)$
for $Z={\cal O}\oplus{\cal O}(1)\oplus{\cal O}(1)$.

Let $Z^{\alpha}$ run over $w^{\mu}$ and $\lambda$, with analogous
definitions for $\hat{U}$. The patching for $TZ\rightarrow Z$ is
\begin{equation}
\hat{V}^{\alpha}=\frac{\partial\hat{Z}^{\alpha}}{\partial Z^{\beta}}V^{\beta}\quad\Rightarrow\quad\begin{array}{c}
\hat{\beta}^{T}=\beta^{T}\\
\hat{\beta}^{A}=\lambda^{-1}\beta^{A}-\lambda^{-2}w^{A}\beta^{\lambda}\\
\hat{\beta}^{\lambda}=-\lambda^{-2}\beta^{\lambda}
\end{array}\,\,,\label{eq:TY0patching}
\end{equation}
and there is one global function of weight zero to consider, $\hat{T}=T$.
The general global section of (\ref{eq:TY0patching}) is
\begin{multline}
\beta=a(T)\frac{\partial}{\partial T}+\left(h^{A}(T)+g^{A}(T)\lambda+j_{\,B}^{A}(T)w^{B}+d(T)\lambda w^{A}+f_{B}(T)w^{B}w^{A}\right)\frac{\partial}{\partial w^{A}}\\
+\left(b(T)+c(T)\lambda+d(T)\lambda^{2}+e_{A}(T)w^{A}+f_{A}(T)w^{A}\lambda\right)\frac{\partial}{\partial\lambda}\,\,,\label{GlobalTY0}
\end{multline}
depending on sixteen (holomorphic) functions of $T$. The global vector
algebra has the decomposition
\[
\check{H}^{0}\left(Z,TZ\right)=\left\{ \underset{\mbox{fifteen}}{\mathfrak{sl}(4,\mathbb{C})}\oplus\underset{\mbox{one}}{\left\{ \frac{\partial}{\partial T}\right\} }\right\} \otimes H\left({\cal O_{T}}\right)
\]
into conformal symmetries of the flat degenerate metric and time translations,
where $H\left({\cal O_{T}}\right)$ are the holomorphic functions
on the time axis.

\subsection{Alpha-Surfaces in the Relativistic Theory}

Via theorem \ref{thm:O(2n)} one can construct some five-dimensional
Riemannian manifolds $M$ as the Kodaira families of curves with normal
bundle ${\cal O}(4)$ in a complex two-fold $Z$, the flat model being
$Z={\cal O}(4)$. In this tangential subsection we'll consider how
the conformal structure on $M$ arises via the presence of alpha surfaces.
The construction is considerably more complicated than that of the
frame advocated throughout this paper, but its existence is nonetheless
reassuring.

The (inhomogeneous) patching for ${\cal O}(4)$ is
\[
\hat{S}=\lambda^{-4}S
\]
and so the global sections are
\[
S|=t+4u\lambda+6x\lambda^{2}+4v\lambda^{3}+w\lambda^{4}
\]
for $x^{a}=(t,u,x,v,w)\in M$.

The direct calculation of the frame (the case $n=2$ in theorem \ref{thm:O(2n)})
is trivial in the flat case and clearly gives us $e^{A^{\prime}B^{\prime}C^{\prime}D^{\prime}}=dx^{A^{\prime}B^{\prime}C^{\prime}D^{\prime}}$
and hence every tensor in the span of the frame, but here we'll consider
alpha surfaces.

The null vectors $\delta x^{a}$ are those for which
\[
0=\delta t+4\delta u\lambda+6\delta x\lambda^{2}+4\delta v\lambda^{3}+\delta w\lambda^{4}
\]
has a unique solution in $\lambda$. The classical theory of quartic
equations tells us that this happens when $\delta x^{a}$ solves simultaneously
the three conditions
\begin{multline*}
\Delta_{6}=-\delta t^{3}\delta w^{3}+12\delta t^{2}\delta u\delta v\delta w^{2}+27\delta t^{2}\delta v^{4}-54\delta t^{2}\delta v^{2}\delta w\delta x\\
+18\delta t^{2}\delta w^{2}\delta x^{2}+6\delta t\delta u^{2}\delta v^{2}\delta w-54\delta t\delta u^{2}\delta w^{2}\delta x-108\delta t\delta u\delta v^{3}\delta x\\
+180\delta t\delta u\delta v\delta w\delta x^{2}+54\delta t\delta v^{2}\delta x^{3}-81\delta t\delta w\delta x^{4}+27\delta u^{4}\delta w^{2}\\
+64\delta u^{3}\delta v^{3}-108\delta u^{3}\delta v\delta w\delta x-36\delta u^{2}\delta v^{2}\delta x^{2}+54\delta u^{2}\delta w\delta x^{3}=0
\end{multline*}
\[
\Delta_{4}=-\delta t\delta w^{3}+4\delta u\delta v\delta w^{2}+12\delta v^{4}-24\delta v^{2}\delta w\delta x+9\delta w^{2}\delta x^{2}=0
\]
\[
\Delta_{2}=-\delta t\delta w+4\delta u\delta v-3\delta x^{2}=0.
\]
\textbf{Claim}: The vanishing of $\Delta_{6}$, $\Delta_{4}$, and
$\Delta_{2}$ is equivalent to $\delta x^{a}$ falling into the union
of the kernels of the span of the frame $e^{A^{\prime}B^{\prime}C^{\prime}D^{\prime}}$.
\begin{itemize}
\item The latter condition $\Delta_{2}=0$ is exactly what one would expect
for the metric, giving us the symmetric two-form $g$ above. (It agrees
exactly with the $g$ one would calculate from the direct frame method.)
Concretely,
\[
[g]=-\delta t\delta w+4\delta u\delta v-3\delta x^{2}.
\]
 
\item The vanishing of $\Delta_{2}$ and $\Delta_{6}$ simultaneously is
equivalent to $\delta x^{a}$ lying in the kernel of both $g$ and
a symmetric three-form $\mathcal{G}_{3}$ (which is also consistent
with the direct frame calculation). In particular
\[
\left[{\cal G}_{3}\left(\delta x^{a},\delta x^{a},\delta x^{a}\right)\right]^{2}\propto\Delta_{6}\qquad\mbox{when}\quad\mbox{\ensuremath{\Delta_{2}}=0}\,,
\]
for
\[
{\cal G}_{3}\left(\delta x^{a},\delta x^{a},\delta x^{a}\right)=-\delta t\delta v^{2}+\delta t\delta w\delta x-\delta u^{2}\delta w+2\delta u\delta v\delta x-\delta x^{3}.
\]
\item Finally, when $\Delta_{2}=\Delta_{6}=0$ the vanishing of $\Delta_{4}$
is equivalent to having
\begin{equation}
\left[\delta x\delta w-\delta v^{2}=0\quad\mbox{and}\quad\delta t\delta x-\delta u^{2}=0\right]\quad\mbox{and either}\quad\delta w\delta t-\delta x^{2}=0\quad\mbox{or}\quad\delta u\delta v-\delta x^{2}=0.\label{degensO4}
\end{equation}
These (effectively three) conditions are the requirement that $\delta x^{a}$
lie in the kernel of three rank-three symmetric two-forms which we
identify as $e_{\qquad A^{\prime}B^{\prime}}^{0^{\prime}0^{\prime}}\otimes e^{0^{\prime}0^{\prime}A^{\prime}B^{\prime}}$,
$e_{\qquad A^{\prime}B^{\prime}}^{0^{\prime}1^{\prime}}\otimes e^{0^{\prime}1^{\prime}A^{\prime}B^{\prime}}$,
and $e_{\qquad A^{\prime}B^{\prime}}^{1^{\prime}1^{\prime}}\otimes e^{1^{\prime}1^{\prime}A^{\prime}B^{\prime}}$.
\item The rest of the canonical symmetric two-forms $e_{\qquad A^{\prime}B^{\prime}}^{C^{\prime}D^{\prime}}\otimes e^{E^{\prime}F^{\prime}A^{\prime}B^{\prime}}$
also arise, but as redundant conditions equivalent to (\ref{degensO4}).
(One could choose to isolate three other conditions, say one rank-three
symmetric two-form and two rank-four symmetric two-forms, and fit
those to canonical forms instead, but for concreteness we have chosen
the three rank-three symmetric two-forms.)
\end{itemize}
Thus we can obtain the induced metric via either direct calculation
or (in a more complicated fashion) by the usual twistor theory arguments.
Once one has the frame one has every canonical form discussed above.

\section*{Acknowledgments}

I would like to thank Maciej Dunajski for stimulating discussions
and I am grateful to STFC and DAMTP for financial support.

\end{document}